\documentclass[oneside,english]{amsart}
\usepackage[T1]{fontenc}
\usepackage[latin9]{inputenc}
\usepackage{mathtools}
\usepackage{enumitem}
\usepackage{amstext}
\usepackage{amsthm}
\usepackage{amssymb}

\usepackage{hyperref}
\usepackage{caption}

\usepackage{verbatim}
\usepackage{units}

\makeatletter
\numberwithin{equation}{section}
\theoremstyle{plain}
\newtheorem{thm}{\protect\theoremname}
\theoremstyle{definition}
\newtheorem{defn}[thm]{\protect\definitionname}
\theoremstyle{remark}
\newtheorem*{rem*}{\protect\remarkname}
\theoremstyle{plain}

\theoremstyle{remark}
\newtheorem{rem}[thm]{\protect\remarkname}
\theoremstyle{plain}
\newtheorem{lem}[thm]{\protect\lemmaname}
\newlist{casenv}{enumerate}{4}
\setlist[casenv]{leftmargin=*,align=left,widest={iiii}}
\setlist[casenv,1]{label={{\itshape\ \casename} \arabic*.},ref=\arabic*}
\setlist[casenv,2]{label={{\itshape\ \casename} \roman*.},ref=\roman*}
\setlist[casenv,3]{label={{\itshape\ \casename\ \alph*.}},ref=\alph*}
\setlist[casenv,4]{label={{\itshape\ \casename} \arabic*.},ref=\arabic*}
\theoremstyle{plain}
\newtheorem{cor}[thm]{\protect\corollaryname}
\theoremstyle{plain}
\newtheorem{prop}[thm]{\protect\propositionname}
\theoremstyle{plain}
\newtheorem*{cor*}{\protect\corollaryname}
\theoremstyle{plain}
\newtheorem*{lem*}{\protect\lemmaname}

\newtheorem{quest}{Question}

\def\R{\mathbb R}

\def\T{\mathbb T}

\newcommand{\trees}{\mathcal{T}\kern-.5mm rees}

\usepackage{babel}
\providecommand{\lemmaname}{Lemma}
\providecommand{\propositionname}{Proposition}
\providecommand{\theoremname}{Theorem}

\makeatother

\usepackage{babel}
\providecommand{\casename}{Case}
\providecommand{\corollaryname}{Corollary}
\providecommand{\definitionname}{Definition}
\providecommand{\factname}{Fact}
\providecommand{\claimname}{Claim}
\providecommand{\lemmaname}{Lemma}
\providecommand{\propositionname}{Proposition}
\providecommand{\remarkname}{Remark}
\providecommand{\theoremname}{Theorem}

\begin{document}

\title{Non-classifiability of mixing zero-entropy diffeomorphisms up to isomorphism}
\author{Marlies Gerber$^1$}   
	\thanks{$^1$ Indiana University, Department of Mathematics, Bloomington, IN 47405, USA}
	\author{Philipp Kunde$^2$} 
	\thanks{$^2$ Oregon State University, Department of Mathematics, Corvallis, OR 97331, USA.}

    \begin{abstract}
We show that the problem of classifying, up to isomorphism, the collection of zero-entropy mixing automorphisms of a standard non-atomic probability space, is intractible. More precisely, the collection of isomorphic pairs of automorphisms in this class is not Borel, when considered as a subset of the Cartesian product of the collection of measure-preserving automorphisms with itself.  
This remains true if we restrict to zero-entropy mixing automorphisms that are also $C^{\infty}$ diffeomorphisms of the five-dimensional torus. In addition, both of these results still hold if ``isomorphism'' is replaced by ``Kakutani equivalence.'' 

In our argument we show that for a uniquely and totally ergodic automorphism $U$ and a particular family of automorphisms $\mathcal{S}$, if $T\times U$ is isomorphic to $T^{-1}\times U$ with $T\in\mathcal{S}$ then $T$ is isomorphic to ${T^{-1}}$. However, this type of ``cancellation'' of factors from isomorphic Cartesian products is not true in general. We present an example due to M. Lema\'nczyk of two weakly mixing automorphisms $T$ and $S$ and an irrational rotation $R$ such that $T\times R$ is isomorphic to $S\times R$, but $T$ and $S$ are not isomorphic.

\end{abstract}
        
\maketitle

\insert\footins{\footnotesize - \\
\textit{2020 Mathematics Subject classification:} Primary: 37A35; Secondary: 37A20, 37A05, 37A25, 37C40, 03E15\\
\textit{Key words: } Mixing, measure-theoretic isomorphism, Kakutani equivalence, anti-classification, complete analytic, smooth ergodic theory}

\section{\label{sec:Introduction}Introduction}

We consider the problem, first posed by J. von Neumann \cite{N32}, of classifying ergodic automorphisms up to isomorphism.
Let $(X,\mathcal{B},\mu)$ and $(Y,\mathcal{C},\nu)$ be standard non-atomic probability spaces. If $T$ and $S$ are measure-preserving automorphisms of $X$ and $Y$, respectively, then we say that $T$ is {\it isomorphic} to $S$, written 
$T \cong S$, if there exists a measure-preserving bijection $\varphi:X\to Y$ such that $\varphi\circ T=S\circ\varphi$ (modulo sets of measure zero). 

Let $\mathcal{X}$ denote the collection of measure-preserving automorphisms of a particular standard non-atomic probability space $(X,\mathcal{B},\mu).$ We endow  $\mathcal{X}$ with the weak topology, that is, for $T_n, T$ in $\mathcal{X}$, $T_n\to T$ if and only if 
$$\mu(T_nA\Delta TA)\to 0 {\text{ for every}} A\in\mathcal{B}.$$
Then $\mathcal{X}$ is a Polish space, which means that the topology is compatible with a complete separable metric. We let $\mathcal{E}$ denote the collection of ergodic elements of $\mathcal{X}$, which is a dense $G_{\delta}$ subset of $\mathcal{X}$. 

The most important and well-known classification results are the Halmos-von Neumann \cite{HN42}  classification of ergodic automorphisms with discrete spectrum and Ornstein's classification of Bernoulli shifts \cite{O70}. However, for isomorphism of ergodic automorphisms in general, no such classification results exist. In fact,  M. Foreman, D. Rudolph, and B. Weiss \cite{FRW11} showed that ergodic measure-preserving automorphisms of $[0,1]$ {\it cannot} be classified up to isomorphism. More precisely, they showed that the collection $\mathcal{R}:=\{(T,S)\in \mathcal{E}\times \mathcal{E}: T\cong S\}$ is a complete analytic subset of $\mathcal{X}\times\mathcal{X}$, and therefore $\mathcal{R}$ is not Borel in $\mathcal{X}\times\mathcal{X}$. Thus there is no method starting with basic open sets in $\mathcal{X}\times\mathcal{X}$ that can reliably determine, in countably many steps, whether or not two elements $T,S\in\mathcal{E}$ are isomorphic. An earlier work of G. Hjorth \cite{H01} proved the same result for $\mathcal{E}$ replaced by $\mathcal{X}$. However, Hjorth's work used nonergodic transformations in an essential way. Prior to \cite{FRW11}, it seemed that the complexity of the isomorphism equivalence relation, as demonstrated in \cite{H01}, might be due to the complexity of the way in which the ergodic components are assembled, rather than the complexity of the ergodic transformations themselves. 

For the purposes of this introduction, we will say that a Borel collection of automorphisms (or diffeomorphisms) is {\it non-classifiable} up to an equivalence relation $\mathcal{R}$ if the set of $\mathcal{R}$-equivalent pairs in this collection is not Borel. This is a strong form of non-classifiability. An example of a weaker form of non-classifiability is the absence of any complete numerical Borel invariant. As was shown by Foreman and A. Louveau in 1995, the Halmos-von Neumann classification already implies that there is no complete numerical Borel invariant for discrete spectrum transformations \cite[Sect. 5]{Fsurvey}, although the spectrum gives a complete algebraic Borel invariant. An earlier proof of the absence of a complete numerical Borel invariant for a family of transformations is J. Feldman's observation \cite{Fe74} that no such invariant exists for the uncountable family of Kolmogorov automorphisms (or K-automorphisms for brevity) that are not Bernoulli constructed by D. Ornstein and P. Shields \cite{OS73}. 

Even though ergodic automorphisms are non-classifiable up to isomorphism, we can ask whether there are any natural Borel subcollections that are classifiable. As we discuss below, weakly mixing automorphisms, K-automorphisms, and zero-entropy mixing automorphisms are all non-classifiable up to isomorphism. On the positive side, Foreman, Rudolph, and Weiss \cite{FRW11} used a theorem of J. King \cite{Ki86} to prove that the rank-one automorphisms, which form a dense $G$-delta subset of $\mathcal{E}$, are classifiable up to isomorphism, but until recently there was no useful structure theorem to classify them. Such a structure theorem has now been obtained by A. Danilenko and M. Vieprik \cite{DVpp}.

The smooth realization problem, described by A. Katok as one of the five most resistant problems in dynamics \cite{Kpp}, is whether every ergodic element of $\mathcal{X}$ with finite entropy can be realized as (is isomorphic to) a smooth diffeomorphism of a compact manifold with a smooth invariant measure. In this context, it is interesting to know whether non-classifiability results for isomorphism of ergodic automorphisms of finite entropy have analogs for diffeomorphisms of compact manifolds. If, for instance, it could be shown that a Borel family consisting of those finite entropy ergodic automorphisms with some specific additional properties invariant under isomorphism, such as the zero-entropy and mixing properties, are not classifiable, but smooth ergodic measure-preserving diffeomorphisms with these properties can be classified, then we would have a negative answer to the smooth realization problem: There would not exist a Borel map from finite-entropy automorphisms to smooth realizations of these automorphisms. This strategy could also be applied to other equivalence relations. As we will see in the discussion below, this has not yet led to any such negative result, which gives some weak evidence that the answer to the smooth realization problem might be positive. 

Foreman and Weiss \cite{FW22} showed that smooth ergodic diffeomorphisms of the two-dimensional disk, annulus, or torus are not classifiable up to isomorphism. Thus, in the case of general ergodic automorphisms we cannot use the above strategy to obtain a negative answer to the smooth realization problem. Note that the weak topology on smooth diffeomorphisms of a compact manifold that preserve a smooth measure is weaker than the Diff$^{\infty}$ topology. Therefore non-classifiability in the smooth setting implies non-classifiability in the automorphism setting. 

Another well-known equivalence relation for ergodic measure-preserving automorphisms is Kakutani equivalence, which is a weaker equivalence relation than isomorphism. Two ergodic measure-preserving automorphisms $T$ and $S$ of a probability space $X$ are said to be {\it Kakutani equivalent}, written $T\sim S$, if there exist positive measure subsets $A$ and $B$ of $X$ such that the first return maps $T_A$ of $T$ to $A$, and $S_B$ of $S$ to $B$, are isomorphic with respect to the normalized restrictions of the probability measure to $A$ and $B$, respectively. Or equivalently, $T$ and $S$ are isomorphic to cross-sections of the same ergodic flow \cite{Ka43}. In \cite{GK25}, we showed that ergodic automorphisms are non-classifiable up to Kakutani equivalence, and this is also true in the smooth setting. Our approach, which is partially based on that in \cite{FRW11}, also gives these results for isomorphism. The second author extended the results of \cite{FRW11} and \cite{FW22} to the weakly mixing setting, both for isomorphism and Kakutani equivalence \cite{Ku24}. Moreover, we showed in \cite{GKpp} that  K-diffeomorphisms (that is, K-automorphisms that are also diffeomorphisms preserving a smooth measure on a compact manifold) are non-classifiable up to isomorphism and Kakutani equivalence. 
Even though K-automorphisms are mixing of all orders, they must have positive entropy. The zero-entropy mixing case that we treat in the present paper utilizes new methods, namely B. Fayad's construction of mixing flows and a certain cancellation property of factors in Cartesian products, that were not previously used in non-classifiability results. 

The non-classifiability of ergodic automorphisms and  diffeomorphisms up to Kakutani equivalence obtained in \cite{GK25} is a key ingredient in obtaining non-classifiability of zero-entropy mixing automorphisms and diffeomorphisms up to isomorphism as well as Kakutani equivalence. In particular, we use two families of transformations $\{F_i(\mathcal{T}):\mathcal{T}\in\mathcal{T} rees\}$, $i=1,2$, from \cite{GK25}, which are described in Section \ref{sec:Prelim}. The family $\{F_1(\mathcal{T})\}$ is used for the automorphism case and $\{F_2(\mathcal{T})\}$ is used for the diffeomorphism case. 

Fayad constructed mixing flows that are special flows under a function and over a base consisting of Cartesian products of the form $U\times R_{\alpha,\alpha'}$, provided that $U\times R_{\alpha,\alpha'}$ 
is uniquely ergodic \cite[Theorem 3]{Fa02}. Here $R_{\alpha,\alpha'}$ is a rotation by $(\alpha,\alpha')$ on $\mathbb{T}^2$ for carefully chosen rationally independent $\alpha,\alpha'$. In fact, by an argument due to A. Kanigowski that is presented in our Proposition \ref{prop:Fayad generalization}, Fayad's theorem still holds even if $U\times R_{\alpha,\alpha'}$ is just assumed to be ergodic, not necessarily uniquely ergodic. 

Once we obtain mixing flows from this construction, the most difficult remaining step (done in Proposition \ref{prop:cancel} and Theorem \ref{thm:main}) for showing non-classifiability of zero-entropy mixing automorphisms and diffeomorphisms is to show that for $i=1,2$ and the ergodic skew product $U$ over $R_{\alpha,\alpha'}$ that occurs in Theorem \ref{thm:main}, if $T=F_i(\mathcal{T})$ for $\mathcal{T}\in\mathcal{T} rees$ and $T\times U\sim T^{-1}\times U$, then $T\cong T^{-1}$. 

In a forthcoming work, P.~Kucharski shows that Fayad's flows are actually mixing of all orders. This would allow the generalization of our anti-classification results to zero-entropy systems that are mixing of all orders.

In Corollary \ref{cor:general n}, which is not used in the proof of our main theorems,
we show how Proposition \ref{prop:cancel} implies that ergodic diffeomorphisms of $\T^3$ are non-classifiable up to isomorphism or up to Kakutani equivalence. This is not surprising, given the results for $\T^2$ \cite[GK25]{FW22}, but it is interesting that the $\T^3$ case can be handled without very much additional work. Non-classifiability for $\mathbb{T}^n$ with $n\ge 4$ can also be obtained from Corollary \ref{cor:general n}. Alternatively, the isomorphism version follows from the $\T^2$ and $\T^3$ cases together with the results of J.-P. Thouvenot and T. Austin that are described below. 

Our Proposition \ref{prop:cancel} answers a special case of the following general question, which is discussed in Sect. \ref{sec:Cancel}.

\begin{quest}
\label{Question}
Suppose $T$ and $S$ are automorphisms such that $T\times U\cong S\times U$ for a particular automorphism $U$, where $T\times U$ is ergodic. Does it follow that $T\cong S$? What if $\cong$ is replaced by $\sim$ in the hypothesis and the conclusion? 
\end{quest}

For general $U$, the answer to both of these questions is ``no''. In Sect. \ref{sec:Cancel} we give an example of $T,S,U$ such that $T\times U\cong S\times U$, but $T\nsim S.$ On the other hand, it follows from work of Thouvenot \cite{T08} and Austin \cite{Au18} that for $B$ equal to a Bernoulli shift and $S$ and $T$ finite entropy transformations, $T\times B\cong S\times B$ implies $T\cong S$. If we let $B$ be a finite-entropy Bernoulli shift that is also a diffeomorphism of $\mathbb{T}^2$ preserving Lebesgue measure $\lambda_2$, this result shows that for any manifold $M$ with a smooth measure $\mu$, classification up to isomorphism of diffeomorphisms of $M$ that preserve $\lambda_2$ can be reduced to the classification up to isomorphism of diffeomorphisms of $M\times\mathbb{T}^2$ that preserve $\mu\times\lambda_2$. In other words, the classification problem is at least as complex for $M\times \mathbb{T}^2$ as it is for $M$, as expected.

If we would like to apply a similar strategy for $M\times S^1$, then we might attempt to show that for an irrational rotation $R$ on $S^1$, $T\times R\cong S\times R$ implies $T\cong S$, at least when $T\times R$ is ergodic. However, this is not true in general. In Sect. 5, we present an example due to M. Lema\'nczyk of two weakly mixing transformations $T_{\varphi}$ and $T_{\psi}$ and an irrational rotation $R$ of $S^1$ such that $T_{\varphi}\times R\cong T_{\psi}\times R$, but $T_{\varphi}\ncong T_{\psi}$. 
If the rotation on $S^1$ is replaced with the flip $f:\{0,1\}\to \{0,1\}$ by $f(0)=1, f(1)=0$, and $T$ is Rudolph's example \cite{Ru76} of a K-automorphism that has two non-isomorphic square roots $S_1$ and $S_2$, then $S_1\times f$ and $S_2\times f$ are isomorphic, thus providing a negative answer to the isomorphism version of the question when the circle rotation is replaced by $f$. 

In an unpublished manuscript, A. Rosenthal \cite{Ro82} constructed two $K$-flows that are isomorphic at time $1$ and not isomorphic at time $\alpha,$ for some irrational $\alpha$. If $T$ and $S$ are the $K$-automorphisms that constitute the time $\alpha$ maps of these flows and $R$ is the rotation of a circle by angle $2\pi\alpha$, then it is easy to verify (see \cite[p.2]{Ro82}) that $T\times R$ is isomorphic to $S\times R.$ Rosenthal used ideas of Rudolph \cite{Ru76} and Ornstein-Shields \cite{OS73} in the construction of these $K$-flows.

The converse of the question above for Kakutani equivalence has a negative answer, because it is possible to have automorphisms $T$ and $S$, both of entropy 0, and a finite entropy Bernoulli shift $B$ such that $T\sim S$, but $T\times B\nsim S\times B$ \cite[p.101]{ORW82}. The following related question stated in \cite{ORW82} still seems to be open: Suppose $T$ and $S$ are automorphisms of entropy 0 and $B$ is a finite entropy Bernoulli shift such that $T\times B \sim S\times B$. Is $T\sim S$?

\section{\label{sec:Statements}Statements of Main Results}
As in the introduction, we let $\mathcal{X}$ be the set of measure-preserving automorphisms of a particular standard non-atomic probability space $(Y,\mu)$ and we endow $\mathcal{X}$ with the weak topology.
The collection of zero-entropy mixing automorphisms of $(Y,\mu)$, which we denote by $\mathcal{ZM}$, is a Borel subset of $\mathcal{X}.$ However,  we show that the isomorphism equivalence relation, restricted to $\mathcal{ZM}\times\mathcal{ZM}$ is a complete analytic set, and hence is not Borel. 

We use $\lambda_k$ (or sometimes just $\lambda$) to denote $k$-dimensional Lebesgue measure on the $k$-dimensional torus $\mathbb{T}^k$, the $k$-dimensional product of the unit interval $[0,1]^k$, or $\mathbb{R}^k$. In the smooth setting we consider those zero-entropy mixing automorphisms that are also $C^{\infty}$ diffeomorphisms of $\mathbb{T}^5$ that preserve $\lambda_5.$ The collection of such automorphisms will be denoted by $\mathcal{ZM}\ \cap$ Diff$^{\infty}(\T^5 \mathbb,\ \lambda_5),$ and this collection is endowed with the Diff$^{\infty}$ topology. 

Note that the weak topology is weaker than the Diff$^{\infty}$ topology on $\mathcal{ZM}\text{\ }\cap $ 

\noindent Diff$^{\infty}(\mathbb{T}^5,\lambda_5).$ Thus sets that are not Borel in the Diff$^{\infty}$ topology are also not Borel in the weak topology. In particular, the non-Borelness in the conclusion of Theorem \ref{Main Automorphism} follows from the non-Borelness in the conclusion of Theorem \ref{Main Smooth}. However, the proof of Theorem \ref{Main Automorphism} by itself is easier than the proof of Theorem \ref{Main Smooth}. The reader who is only interested in non-classifiability of zero-entropy mixing {\it automorphisms} may omit Propositions \ref{prop:Fayad generalization} and \ref{prop:isotopy}. In the automorphism case, the proof of Proposition \ref{prop:ue} is needed only in the case $i=1,$ which is simpler than the case $i=2$, and the proof actually shows the unique ergodicity of $F_1(\mathcal{T})\times R_{\alpha,\alpha'}.$  

We first state our results in the automorphism case. 
\begin{thm}
\label{Main Automorphism}
The sets $$\mathcal{S}_{\text{is}}:=\{(S,T):S,T\in \mathcal{ZM} \text{\ and\ } S\cong T\}$$
and $$\mathcal{S}_{\text{Ka}}:=\{(S,T):S,T\in \mathcal{ZM} \text{\ and\ } S\sim T\}$$
are complete analytic subsets of $\mathcal{ZM}\times\mathcal{ZM}.$ In particular they are not Borel.    
\end{thm}
In the smooth setting we have the following result.
\begin{thm}
    \label{Main Smooth}
The sets 
$$\mathcal{S}^{\infty}_{\text{is}}:=\{(S,T):S,T\in \mathcal{ZM}\cap\text{Diff}^{\infty}(\mathbb{T}^5,\lambda) \text{\ and\ }S\cong T\}$$ and 
$$\mathcal{S}^{\infty}_{\text{Ka}}:=\{(S,T):S,T\in \mathcal{ZM}\cap\text{Diff}^{\infty}(\mathbb{T}^5,\lambda) \text{\ and\ }S\sim T\}$$ 
are complete analytic subsets of $(\mathcal{ZM}\cap\text{Diff}^{\infty}(\mathbb{T}^5,\lambda))
\times (\mathcal{ZM}\cap\text{Diff }^{\infty}(\mathbb{T}^5,\lambda)).$ In particular, they are not Borel.

\end{thm}

Theorems \ref{Main Automorphism} and \ref{Main Smooth} are, respectively, corollaries of the Theorems \ref{AutoMix} and \ref{SmoothMix} below (see Sect. \ref{sec:Prelim}). We will use the Theorems \ref{AutoErgodic} and \ref{SmoothErgodic} from \cite{GK25}, together with results and techniques from \cite{Fa02} that we describe in Section \ref{sec:Proof}. Theorem \ref{AutoErgodic} is Proposition 25 in \cite{GK25}, with $F_1=\Psi,$ and Theorem \ref{SmoothErgodic} is proved in Section 10 of \cite{GK25}, with $F_2 = F^s.$

\begin{thm}
\label{AutoErgodic}
There is a continuous map 
$$F_1:\mathcal{T} rees\to\mathcal{E}$$
such that 
\begin{itemize}
    \item[(i)] If $\mathcal{T}$ has an infinite branch, then $F_1(\mathcal{T})$ is isomorphic (and hence Kakutani equivalent) to $(F_1 (\mathcal{T}))^{-1}.$
    \item[(ii)] If $\mathcal{T}$ does not have an infinite branch, then $F_1(\mathcal{T})$ is not Kakutani equivalent (and hence not isomorphic) to $(F_1 (\mathcal{T}))^{-1}.$
\end{itemize}
\end{thm}

\begin{thm}
\label{SmoothErgodic} 
There is a continuous map
$$F_2:\mathcal{T} rees\to\mathcal{E}\cap \text{Diff}^{\infty}(\mathbb{T}^2,\lambda_2)$$
such that (i) and (ii) of Theorem \ref{AutoErgodic} hold for $F_1$ replaced by $F_2$.

\end{thm}

For the maps $F_i$, $i=1,2,$ constructed in \cite{GK25}, each $F_i(\mathcal{T})$ has entropy zero, but is {\it{not}} weakly mixing. We will use \cite{Fa02} to upgrade Theorems \ref{AutoErgodic} and \ref{SmoothErgodic} to $\mathcal{ZM},$ but with $\mathbb{T}^2$ replaced by $\mathbb{T}^5.$ 

\begin{thm}
\label{AutoMix} 
There is a continuous map
$$\widetilde{F}_1=(\widetilde{F}_{1,a},\widetilde{F}_{1,b}):\mathcal{T} rees\to \mathcal{ZM}\times\mathcal{ZM}$$ such that 

\begin{itemize}
    \item[(i)] If $\mathcal{T}$ has an infinite branch, then $\widetilde{F}_{1,a}(\mathcal{T})$ is isomorphic (and hence Kakutani equivalent) to $\widetilde{F}_{1,b} (\mathcal{T}).$
    \item[(ii)] If $\mathcal{T}$ does not have an infinite branch, then $\widetilde{F}_{1,a}(\mathcal{T})$ is not Kakutani equivalent (and hence not isomorphic) to $\widetilde{F}_{1,b} (\mathcal{T}).$
   
\end{itemize}
\end{thm}

\begin{thm}
    \label{SmoothMix}
  There is a continuous map
  $$\widetilde{F}_2=(\widetilde{F}_{2,a},\widetilde{F}_{2,b}):\mathcal{T} rees\to\big(\mathcal{ZM}\cap\text{Diff}^{\infty}(\mathbb{T}^5,\lambda)\big)\times \big(\mathcal{ZM}\cap\text{Diff}^{\infty}(\mathbb{T}^5,\lambda)\big)$$ such that (i) and (ii) of Theorem \ref{AutoMix} hold with $\widetilde{F}_{1,a}$ and $\widetilde{F}_{1,b}$ replaced by $\widetilde{F}_{2,a}$ and $\widetilde{F}_{2,b},$ respectively.
\end{thm}

\section{\label{sec:Prelim}Preliminaries}
\subsection{\label{subsec:Borel}Borel Reductions}
To prove our anti-classification results we follow the general strategy from \cite{FRW11}. An important tool is the concept of a
reduction. 
\begin{defn}
Let $X$ and $Y$ be Polish spaces and $A\subseteq X$, $B\subseteq Y$.
A function $f:X\to Y$ \emph{reduces} $A$ to $B$ if and only if
for all $x\in X$: $x\in A$ if and only if $f(x)\in B$. 

Such a function $f$ is called a \emph{Borel (respectively, continuous) reduction}
if $f$ is a Borel (respectively, continuous) function.
\end{defn}

$A$ being reducible to $B$ can be interpreted as saying that $B$
is at least as complicated as $A$. We note that if $B$ is Borel
and $f$ is a Borel reduction, then $A$ is also Borel. Equivalently, if $f$ is a Borel reduction of $A$ to $B$ and $A$ is not Borel,
then $B$ is not Borel.

\begin{defn}
If $X$ is a Polish space and $B\subseteq X$, then $B$ is \emph{analytic}
if and only if it is the continuous image of a Borel subset of a Polish
space. Equivalently, there is a Polish space $Y$ and a Borel set
$C\subseteq X\times Y$ such that $B$ is the $X$-projection of $C$.
\end{defn}

\begin{defn}
An analytic subset $A$ of a Polish space $X$ is called \emph{complete analytic} if every analytic set can be continuously reduced to $A$.
\end{defn}
If $A$ is continuously reducible to an analytic set $B$ and $A$ is a complete analytic set, then so is $B$.
There are analytic sets that are not Borel  (see, for example, \cite[Sect. 14]{Ke95}). This implies that a complete analytic set is not Borel. 

In the Polish space $\mathcal{T}rees$, the subset consisting of trees with an infinite branch is an example of a complete analytic set \cite[Sect.~27]{Ke95}.

  We now define trees, the space  $\mathcal{T}\kern-.5mm rees$, and infinite branches of trees.

 Let $\mathbb{N}$ be the nonnegative integers, and let $\mathbb{N}^{<\mathbb{N}}$ denote the collection 
 of finite sequences in $\mathbb{N}$. By convention, the empty sequence $\emptyset$ is in $\mathbb{N}^{<\mathbb{N}}$.

 \begin{defn}
 A \emph{tree} is a set $\mathcal{T}\subseteq\mathbb{N}^{<\mathbb{N}}$
	such that if $\tau=\left(\tau_{1},\dots,\tau_{n}\right)\in\mathcal{T}$
	and $\sigma=\left(\tau_{1},\dots,\tau_{m}\right)$ with $m\leq n$
	is an initial segment of $\tau$, then $\sigma\in\mathcal{T}$.
 Each tree $\mathcal{T}$ contains $\emptyset$, which is regarded as
 an initial segment of all $\tau\in\mathcal{T}$.
 \end{defn}
 
	\begin{defn} 
		A function $f:\{0,\dots,n-1\}\to \mathbb{N}$ such that $(f(0),\dots,f(n-1))\in\mathcal{T}$ is called a \emph{branch} of $\mathcal{T}$ of length $n$. 
         The collection  $\mathcal{T}\kern-.5mm rees$
	       consists of those trees that contain branches of arbitrarily long length.
        An \emph{infinite branch} through $\mathcal{T}$ is a function $f:\mathbb{N}\to\mathbb{N}$ such that for all $n\in\mathbb{N}$ we have $\left(f(0),\dots,f(n-1)\right)\in\mathcal{T}$. If a tree has an infinite branch, it is called \emph{ill-founded}.
		If it does not have an infinite branch, it is called \emph{well-founded}.
	\end{defn}

 To describe a topology on the collection of trees, let $\left\{ \sigma_{n}:n\in\mathbb{N}\right\} $
	be an enumeration of $\mathbb{N}^{<\mathbb{N}}$ with the property
	that every 
    initial segment
    of $\sigma_{n}$ is some $\sigma_{m}$
	for $m \leq n$. Under this enumeration subsets $S\subseteq\mathbb{N}^{<\mathbb{N}}$
	can be identified with characteristic functions $\chi_{S}:\mathbb{N}\to\left\{ 0,1\right\} $.
	Such $\chi_{S}$ can be viewed as 
    members of
	an infinite product space $\left\{ 0,1\right\} ^{\mathbb{N}^{<\mathbb{N}}}$
	homeomorphic to the Cantor space. Here, each function $a:\left\{ \sigma_{m}:m<n\right\} \to\left\{ 0,1\right\} $
	determines a basic open set 
	\[
	\left\langle a\right\rangle =\left\{ \chi:\chi\upharpoonright\left\{ \sigma_{m}:m<n\right\} =a\right\} \subseteq\left\{ 0,1\right\} ^{\mathbb{N}^{<\mathbb{N}}}
	\]
	and the collection of all such $\left\langle a\right\rangle $ forms
	a basis for the topology. In this topology the collection of trees is a closed (hence
	compact) subset of $\left\{ 0,1\right\} ^{\mathbb{N}^{<\mathbb{N}}}$. Moreover, the collection $\mathcal{T}\kern-.5mm rees$ 
    is a dense $\mathcal{G}_{\delta}$ subset. Hence, $\mathcal{T}\kern-.5mm rees$ is a Polish space. 
	
We now show how Borel reductions and $\mathcal{T} rees$ are used to prove our main anti-classification results. 

\begin{proof}[Proof of Theorem \ref{Main Automorphism} from Theorem \ref{AutoMix}]
The map $\widetilde{F}_1$ in Theorem \ref{AutoMix} provides a continuous reduction from the set of trees with an infinite branch to $\mathcal{S}_{\text {is}}=\{(S,T)\in \mathcal{ZM}\times\mathcal{ZM}:S\text{\  is isomorphic to\  } T\}$. As observed in \cite{FRW11}, $\{(S,T)\in \mathcal{X}\times\mathcal{X}: S\text{\  is isomorphic to\  } T\}$ is an analytic set. (For ``isomorphic'' replaced by ``Kakutani equivalent,'' this was shown in \cite[Section 3.2]{GKpp}.) Since $\mathcal{ZM}\times\mathcal{ZM}$ is a Borel subset of the Polish space $\mathcal{X}\times\mathcal{X}$, $\mathcal{S}_{\text {is}}$ is the intersection of a Borel set with an analytic set, which is analytic. Therefore $\mathcal{S}_{\text {is}}$ is a complete analytic set. The proof that $\mathcal{S}_{\text {Ka}}$ is a complete analytic set is proved the same way, using the Kakutani equivalence version of Theorem \ref{AutoMix}. 
\end{proof}

\begin{proof}[Proof of Theorem \ref{Main Smooth} from Theorem \ref{SmoothMix}] The proof is the same as for Theorem \ref{Main Automorphism}, using $\widetilde{F}_2$ instead of $\widetilde{F}_1$ for the reduction.
\end{proof}

\subsection{\label{subsec:Terminology}Some Terminology from Ergodic Theory} 
As in the introduction, we let $(X,\mathcal{B},\mu)$ be
a \emph{standard non-atomic probability space}, that is, there is a bi-measurable measure-preserving bijection (modulo sets of measure zero) from the measure space $(X,\mathcal{B},\mu)$ to $([0,1],\operatorname{Bor},\lambda)$, where $\operatorname{Bor}$ is the collection of Borel sets, and $\lambda$ is Lebesgue measure on $[0,1]$.

A \emph{partition} $\mathcal{P}$ of $(X,\mathcal{B},\mu)$ is a collection $\mathcal{P}=\left\{ c_{\sigma}\right\} _{\sigma\in\Sigma}$
	of subsets $c_{\sigma}\in\mathcal{B}$ with $\mu(c_{\sigma}\cap c_{\sigma'})=0$
	for all $\sigma\neq\sigma'$ and $\mu\left(\bigcup_{\sigma\in\Sigma}c_{\sigma}\right)=1$,
	where $\Sigma$ is a finite set of indices. For two
partitions $\mathcal{P}$ and $\mathcal{Q}$ we define the join of $\mathcal{P}$ and $\mathcal{Q}$ to be the partition $\mathcal{P}\vee \mathcal{Q}=\left\{ c\cap d:c\in \mathcal{P},d\in \mathcal{Q}\right\}$. Then for a sequence of partitions $\{\mathcal{P}_{n}\}_{n=1}^{\infty}$ we inductively define $\vee^N_{n=1}\mathcal{P}_n,$ and we let
$\vee_{n=1}^{\infty}\mathcal{P}_{n}$ be the smallest $\sigma$-algebra containing $\vee^N_{n=1}\mathcal{P}_n$ for each $N \in \mathbb{Z}^+$.
We say that a sequence of partitions $\{\mathcal{P}_{n}\}_{n=1}^{\infty}$ is
a \emph{generating sequence} if $\vee_{n=1}^{\infty}\mathcal{P}_{n} = \mathcal{B}$.
    
For an automorphism $T$ of $(X,\mathcal{B},\mu)$ and a
	partition $\mathcal{P}=\left\{ c_{\sigma}\right\} _{\sigma\in\Sigma}$
	of $X$ we can define the \emph{$(T,\mathcal{P})$-name} for almost every $x\in X$
	by 
	\[
	(a_n)_{n=-\infty}^{\infty}\in\Sigma^{\mathbb{Z}}\text{ with }T^{i}(x)\in c_{a_{i}}.
	\]
If the partition $\mathcal{P}$ is a \emph{generator}
(that is, $\left\{ T^{n}\mathcal{P}\right\} _{n=-\infty}^{\infty}$ is a generating
sequence),
then $(X,\mu,T)$ is isomorphic to $(\Sigma^{\mathbb{Z}},\nu,sh)$
with left shift $sh:\Sigma^{\mathbb{Z}}\to\Sigma^{\mathbb{Z}},$ defined by
\[
sh((x_n)_{n=-\infty}^{\infty})=(x_{n+1})_{n=-\infty}^{\infty},
\]  
and the measure $\nu$ satisfying $\nu\{(x_n)_{n=-\infty}^{\infty}: x_n=a_n \text{\ for\ } n=j,j+1,...,k\}=\mu(T^{-j}(c_{a_j}))\cap T^{-(j+1)}(c_{a_{j+1}})\cap\cdots\cap T^{-k}(c_{a_k})).$ The pair $(T,\mathcal{P})$ is called a \emph{process}. 

    The \emph{Hamming distance}
	between two strings of symbols $a_{1}\dots a_{n}$ and $b_{1}\dots b_{n}$ from a given alphabet $\Sigma$
	is defined by $\overline{d}\left(a_{1}\dots a_{n},b_{1}\dots b_{n}\right)=\frac{1}{n}|\left\{ i:a_{i}\neq b_{i}\right\} |$.
	It can be extended to infinite strings $w=\dots a_{-2}a_{-1}a_{0}a_{1}a_{2}...$
	and $w^{\prime}=\dots b_{-2}b_{-1}b_{0}b_{1}b_{2}\dots$ by $\overline{d}(w,w^{\prime})=\limsup_{n\to\infty}\overline{d}\left(w_{n},w_{n}^{\prime}\right)$,
	where $w_{n}=a_{-n}a_{-n+1}\dots a_{n-1}a_{n}$ and $w_{n}^{\prime}=b_{-n}b_{-n+1}\dots b_{n-1}b_{n}$
	are the truncated strings.

    In the study of Kakutani equivalence, Feldman introduced a weaker notion of distance, now called $\overline{f}$, as a substitute for the Hamming distance $\overline{d}$ in Ornstein's isomorphism theory. To define it, we call a
		collection $\mathcal{M}$ of pairs of indices $(i_{s},j_{s})$, $s=1,\dots,r$
		such that $1\le i_{1}<i_{2}<\cdots<i_{r}\le n$, $1\le j_{1}<j_{2}<\cdots<j_{r}\le m$
		and $a_{i_{s}}=b_{j_{s}}$ for $s=1,2,\dots,r$ a \emph{match} between two strings $a_{1}a_{2}\dots a_{n}$
		and $b_{1}b_{2}\dots b_{m}$. Then 
		\begin{equation}
			\begin{array}{ll}
				\overline{f}(a_{1}a_{2}\dots a_{n},b_{1}b_{2}\dots b_{m})=\hfill\\
				{\displaystyle 1-\frac{2\sup\{|\mathcal{M}|:\mathcal{M}\text{\ is\ a\ match\ between\ }a_{1}a_{2}\cdots a_{n}\text{\ and\ }b_{1}b_{2}\cdots b_{m}\}}{n+m}.}
			\end{array}\label{eq:cl}
		\end{equation}
		We will refer to $\overline{f}(a_{1}a_{2}\cdots a_{n},b_{1}b_{2}\cdots b_{m})$
		as the ``$\overline{f}$-distance'' between $a_{1}a_{2}\cdots a_{n}$
		and $b_{1}b_{2}\cdots b_{m},$ even though $\overline{f}$ does not
		satisfy the triangle inequality unless the strings are all of the
		same length. A match $\mathcal{M}$ is called a\emph{ best possible
			match} if it realizes the supremum in the definition of $\overline{f}$.

        A generalization of Proposition 3.2 in \cite[p.~92]{ORW82} given in \cite[section 9.2]{GK25} guarantees the existence of a consistent sequence of finite codes, that are well-approximating in $\overline{f},$ between two evenly equivalent transformations. Hence, a key step in our later arguments to show that two transformations are not Kakutani equivalent will be to exclude the existence of such a consistent sequence of well-approximating finite codes. Here, a \emph{code} of length $2K+1$ is a function $\phi:\Sigma^{\mathbb{Z}\cap[-K,K]}\to\tilde{\Sigma}$, where $\Sigma,\tilde{\Sigma}$ are alphabets. 
Given a code $\phi$ of length $2K+1,$ the \emph{stationary code}
$\bar{\phi}$ on $\Sigma^{\mathbb{Z}}$ determined by $\phi$ is defined by
\[
    \bar{\phi}(s)(l)=\phi\left(s\upharpoonright[l-K,l+K]\right), \text{ for }s\in\Sigma^{\mathbb{Z}}\text{ and }l\in\mathbb{Z}.
\]

In the following subsection we review the construction of the transformations in \cite{GK25} as symbolic systems. The words are built using specific patterns of blocks called \emph{Feldman patterns} since they originate from Feldman's first example of an ergodic zero-entropy transformation that is not loosely Bernoulli \cite{Fe74}. For $T,N,M\in\mathbb{Z}^{+}$ we define a $(T,N,M)$-Feldman pattern in building
blocks $A_{1},\dots,A_{N}$ of equal length $L$ to be one of the strings
$B_{j}$, $j=1,\dots , M$, given by
\begin{eqnarray}\label{eq:Feldman}
	B_{j}\ \ =   & \left(A_{1}^{TN^{2j}}A_{2}^{TN^{2j}}\dots A_{N}^{TN^{2j}}\right)^{N^{2(M+1-j)}},
\end{eqnarray}
where an exponent on a string indicates the number of times the string is repeated.  Every block $B_{j}$ has total length $TN^{2M+3}L$.
	Moreover, we notice that $B_{j}$ is built with $N^{2(M+1-j)}$ many
	so-called \emph{cycles}: Each cycle winds through all the $N$ building
	blocks.

For $j>i$ the length of $TN^{2j}$ consecutive occurrences of a
building block in $B_j$ is equal to the length of $N^{2(j-i)-1}$ complete cycles in $B_i$. If we assume that distinct building blocks cannot be matched well in $\overline{f}$, this different repetition behavior implies that different Feldman patterns
cannot be matched well in $\overline{f}$ even after a finite code is applied to one of them. This still holds if we just consider   substantial substrings of the two patterns. (We refer to Proposition~42 and Lemma~43 in \cite{GK25} for precise statements. Also see Figure~\ref{fig:fig2} in Subsect. \ref{Main proofs}.) In fact, for this argument, one can replace a cycle in the string to which the code is being applied with any permutation of its building blocks. 

\subsection{\label{subsec:Review}Review of Construction of Maps $F_1$ and $F_2$ in Theorems \ref{AutoErgodic} and \ref{SmoothErgodic}}
The construction of the maps $F_1$ and $F_2$ are described in detail in \cite{GK25}. Below we will just summarize those properties that are needed for the proofs of Theorems \ref{AutoMix} and \ref{SmoothMix}. 
For each $\mathcal{T}\in\mathcal{T} rees$, $F_1(\mathcal{T})$ is described in terms of $n$-words, $\mathcal{W}_n=\mathcal{W}_n(\mathcal{T})$, $n\in\mathbb{N}$, where $n$-words are concatenated to form $(n+1)$-words. We let $0$-words consist of a single symbol in a finite alphabet $\Sigma.$ The following specifications are satisfied. 
\begin{enumerate}[label=(A\arabic*)]
	\item\label{item:A1}  All words in $\mathcal{W}_{n}$ have the same length $h_{n}$.
	\item\label{item:A2}  There are $f_{n},k_n\in\mathbb{Z}^{+}$ such that every word in $\mathcal{W}_{n+1}$ is built by concatenating $k_n$ words
	in $\mathcal{W}_{n}$ and 
	every word in $\mathcal{W}_{n}$ occurs in each word of $\mathcal{W}_{n+1}$
	exactly $f_{n}$ times. Clearly, we have $k_n=f_n|\mathcal{W}_n|$.
	\item\label{item:A3}  If $w=w_{1}\dots w_{k_n}\in \mathcal{W}_{n+1}$ and $w^{\prime}=w_{1}^{\prime}\dots w_{k_n}^{\prime}\in\mathcal{W}_{n+1}$, where $w_{i},w_{i}^{\prime}\in\mathcal{W}_{n}$, then for any $k\geq\lfloor\frac{k_n}{2}\rfloor$ and $1\leq i \leq k_n-k$, we have $w_{i+1}\dots w_{i+k}\neq w_{1}^{\prime}\dots w_{k}^{\prime}$. 
\end{enumerate}

According to the definitions in \cite[Sect. 3.3]{FRW11}, conditions \ref{item:A2} and \ref{item:A3} imply that  $\{\mathcal{W}_n\}_{n\in \mathbb{N}}$
	is a uniquely readable and strongly uniform construction sequence.

The domain of $T=F_1(\mathcal{T})$ is the closed subset $\mathcal{S}\subset \Sigma^{\mathbb{Z}}$ consisting of those doubly infinite sequences $(x_n)_{n=-\infty}^{\infty}$ such that for any $j\le k$, $(x_n)_{n=j}^k$ occurs somewhere as consecutive symbols in some $w\in\cup_{n=0}^{\infty} \mathcal{W}_n$. We let $T$ be the left shift, $sh,$ on $\mathcal{S}$. Clearly $\mathcal{S}$ is shift-invariant, and by \cite[Lemma 11]{FRW11} it follows from \ref{item:A2} and \ref{item:A3} that there is a unique invariant measure $\mu$ for the system $(\mathcal{S},sh)$. 

The topological dynamical system $(\mathcal{S},sh,\mu)$ is isomorphic to an automorphism of $[0,1]$ built from a cutting and stacking construction, described as follows. We start by dividing $[0,1]$ into $|\Sigma|$ intervals of length $1/|\Sigma|$. We make a one-to-one correspondence between these intervals and symbols in $\Sigma,$ and label each interval with the corresponding symbol. At the $n$th stage of the cutting and stacking process, we have a tower of height $h_n$ with $|\mathcal{W}_n|$ columns, all of equal width. There is a one-to-one correspondence between these columns and words in $\mathcal{W}_n$, and the levels of the columns are labeled according to the symbols in the corresponding word of $\mathcal{W}_n$. Each column is divided into subcolumns of equal width, and these subcolumns are stacked to form the tower at stage $n+1$. At every stage of the construction, the levels of the columns in the tower consist of intervals. The transformation $T$ obtained by mapping from each level to the level above it (if there is one) preserves Lebesgue measure, and it is isomorphic to $(\mathcal{S},sh,\mu).$
The partition of $[0,1]$ according to the labels in $\Sigma$ forms a generating partition for the transformation $T$.

The construction of $F_1(\mathcal{T})$ allows the following flexibility, which is essential to our argument.
\begin{enumerate}
    \item[(A4)]\label{item:A4} Once $\mathcal{W}_n$ has been constructed, the parameters in the construction of $\mathcal{W}_{n+1}$  can be chosen so that whenever an $n$-word appears within an $(n+1)$-word, it is part of a block of $\ell_n$ consecutive appearances of that same $n$-word, where $\ell_n$ can be chosen arbitrarily large.  
\end{enumerate}

The construction of $F_2(\mathcal{T})$ is similar, but we allow ``spacers'', $b,e$, which are not in $\Sigma$. These spacers are inserted between groups of $n$-words in the construction of $(n+1)$-words. The collection of $n$-words is now denoted $\mathcal{W}_n^c$, where $c$ stands for ``circular'', because these types of words were introduced in the construction of the so-called circular systems \cite{FW19}. This construction can be done so that property (A4) still holds with $\mathcal{W}_n$ and $\mathcal{W}_{n+1}$ replaced by $\mathcal{W}_n^c$ and $\mathcal{W}_{n+1}^c$. That is, whenever a word in $\mathcal{W}_n^c$ appears in a word in $\mathcal{W}_{n+1}^c$ it is a part of a block of $\ell_n$ consecutive appearances of that same word in $\mathcal{W}_n^c$, with no newly added spacers in between. (Actually, in \cite{FW19} and \cite{GK25} there were $\ell_n-1$ consecutive appearances of each $n$-word, but we find it more convenient to denote the number of consecutive appearances by $\ell_n$ in the present paper.)

The spacers do not significantly affect our arguments, because the proportion, say $\delta_n$, of the length of each $(n+1)$-word that consists of spacers that are added in the construction of the $(n+1)$-word once we already have the $n$-words, can be made as small as we like, by choosing $\ell_n$ sufficiently large. We will require, among other things, that $\Sigma_{n=1} ^{\infty}\delta_n<\infty.$ 

If we construct the subshift $(\mathcal{S},sh)$ of $((\Sigma\cup\{b,e\})^{\mathbb{Z}},sh)$ in the same way as for $F_1(\mathcal{T})$, then $(\mathcal{S},sh)$ is not uniquely ergodic. In particular, the points $(\dots,b,b,b,\dots)$ and $(\dots,e,e,e,\dots)$ are fixed points of $(\mathcal{S},sh).$ 

In the case of $F_2(\mathcal{T})$, we consider the cutting and stacking model of the construction, as we did for $F_1(\mathcal{T})$, except we start with an interval of length $\Pi_{n=1}^{\infty}(1+\delta_n)^{-1}$, and we add intervals corresponding to the spacers that are added in going from $n$-words to $(n+1)$-words, so that the total length of the intervals in the towers increases by a factor of $1+\delta_n$. The invariant probability measure that we consider is Lebesgue measure on the intervals in this increasing union of towers. Copying this measure back over to $(\mathcal{S},sh)$ gives the invariant measure that we use for $(\mathcal{S},sh)$. In \cite[Sect. 10]{GK25} we used the machinery from  \cite{FW19} to obtain an ergodic diffeomorphism of $\mathbb{T}^2$ that preserves Lebesgue measure $\lambda_2$ and is isomorphic to $(\mathcal{S},sh).$ It is this diffeomorphism that we take as our $F_2(\mathcal{T}).$

\section{\label{sec:Proof}Proof of Main Results}

In this section we will give the proofs of Theorems \ref{AutoMix} and \ref{SmoothMix}. Then our main results, Theorems \ref{Main Automorphism} and \ref{Main Smooth}, will follow, as was shown in Subsect. \ref{subsec:Borel}.

If $U$ is an automorphism of a probability space $(X,\mu)$ and $\varphi$ is a function from $X$ to $\R^+$, we let $(U,\varphi)$ denote the flow over $U$ and under $\varphi$ on $X^{\varphi}=\{(x,s):0\le s\le \varphi(x)\},$ where the points $(x,\varphi(x))$ and $(U(x),0)$ are identified. For $t\ge 0$ and $0\le s<\varphi(x)$, the time $t$ map sends $(x,s)$ to:
\[
\begin{cases} (x,t+s),& \text{ if
}t+s<\varphi(x),\\
(U^{n}(x), t+s-\Sigma_{j=0}^{n-1}\varphi(U^j(x))& \text{ if }\Sigma_{j=0}^{n-1}\varphi(U^jx)\le t+s<\Sigma_{j=0}^n \varphi(U^jx).
\end{cases}
\]
This flow preserves $\mu^{\varphi}$, which is the restriction of the product measure $\mu\times\lambda_1$ to $X^{\varphi}$. 

\subsection{\label{sec:Mixing}{Obtaining the mixing property}}
By the work of B. Fayad \cite{Fa02}, there exist positive numbers $\alpha$ and $\alpha'$ (with $1,\alpha,$ and $\alpha'$ linearly independent over the rationals, and the rational convergents of $\alpha,\alpha'$ having exponentially growing denominators) and a real analytic function $\varphi:\mathbb{T}^2\to\mathbb{R}^+$ such that for the rotation $R_{\alpha,\alpha'}$ by $(\alpha,\alpha')$ on $\T^2$, the flow $(R_{\alpha,\alpha'},\varphi)$ is mixing. Moreover, Fayad also obtained mixing flows that have Cartesian products of $R_{\alpha,\alpha'}$ with other dynamical systems in the base. The following Proposition \ref{prop:Fayad generalization} is stated in \cite[Theorem 3]{Fa02}, under the additional hypothesis that $S\times R_{\alpha,\alpha'}$ is uniquely ergodic. A.~Kanigowski pointed out that it suffices to assume that $S\times R_{\alpha,\alpha'}$ is ergodic. 

It is important that $\widetilde{\varphi}$ in Proposition \ref{prop:Fayad generalization} only depends on the coordinates of the second factor in $Y\times\mathbb{T}^2$.  This fact implies that $ (S\times R_{\alpha,\alpha'},\widetilde{\varphi})$ is isomorphic to $(T\times R_{\alpha,\alpha'},\widetilde{\varphi})$ if $S$ is isomorphic to $T$. 
 
\begin{prop}\cite[Theorem 3]{Fa02}\label{prop:Fayad generalization}
Suppose $(Y,\nu)$ is a standard non-atomic probability space, and $S:Y\to Y$ is an automorphism that preserves $\nu.$ Let $\varphi:\mathbb{T}^2\to\mathbb{R}^+$ be as in \cite[Theorem 1]{Fa02} and define $\widetilde{\varphi}:Y\times\mathbb{T}^2\to \mathbb{R}^+$ by $\widetilde{\varphi}(y,z,w)=\varphi(z,w)$ for $(y,z,w)\in Y\times\mathbb{T}^1\times\mathbb{T}^1$. If $(S\times R_{\alpha,\alpha'},\nu\times \lambda_2)$ is ergodic, then the flow $(S\times R_{\alpha,\alpha'},{\widetilde{\varphi}})$ is mixing.
\end{prop}

\begin{proof}[Outline of proof] We present Kanigowski's argument for how  Fayad's proof can be altered to achieve this stronger result (assuming $S\times R_{\alpha,\alpha'}$ is ergodic, but not necessarily uniquely ergodic). 

Let $\{T^t\}$ denote the flow $(S\times R_{\alpha,\alpha'},{\widetilde{\varphi}})$. The key ingredient in the proof that $\{T^t\}$ is mixing is still the uniform stretch property (see \cite[Definition 1]{Fa02}) of the Birkhoff sums of $\widetilde{\varphi}$. For $m\in \mathbb{Z}^+,$ let the $m$th Birkhoff sum for $\widetilde{\varphi}$ be $f_m(y,z,w):=\sum_{i=0}^{m-1}\widetilde{\varphi}((S\times R_{\alpha,\alpha'})^i(y,z,w))=\sum_{i=0}^{m-1}\varphi(R_{\alpha,\alpha'}^i(z,w))$. 
As in \cite{Fa02} we assume $\frac{1}{2}\le\varphi\le\frac{3}{2}$ and the average value of $\varphi$ is 1. 

According to \cite[Lemma 2 and Propositions 4 and 5]{Fa02}, for $t>0$, there exist partial partitions $\eta_t=\{C_i^{(t)}\}$ of $\mathbb{T}^1 $ into intervals satisfying the hypothesis of \cite[Proposition 1]{Fa02} (when adapted to having base transformation $S\times R_{\alpha,\alpha'}$ instead of just $R_{\alpha,\alpha'}$) . That is, 
$$ \sup_{C_i^{(t)}\in\eta_t}|C_i^{(t)}|\underset{t\to\infty}{\longrightarrow} 0$$

$$\sum_{C_i^{(t)}\in\eta_t}|C_i^{(t)}|\underset{t\to\infty}{\longrightarrow}1,$$
and there exist positive functions $\varepsilon(t)$ and $k(t)$ such that
$$\varepsilon(t)\underset{t\to\infty}{\longrightarrow}0 $$

$$k(t)\underset{t\to\infty}{\longrightarrow}\infty,$$ and the function $f=\tilde{\varphi}$ is such that, for any $t$ at least one of the following two conditions is true, where points in $\mathbb{T}^2$ are written as $(z,w)$:
\begin{enumerate}
    \item[(i)] for any $m\in [\frac{1}{2}t,2t],$ for any $w\in\mathbb{T}^1,$ any $y\in Y$ and any $C_i^{(t)},$
    $$f_m(y,\cdot,w)\text{ is }(\varepsilon(t),k(t))\text{-uniformly stretching on }\{y\}\times C_i^{(t)}\times\{w\}; $$

    \item[(ii)]for any $m\in [\frac{1}{2}t,2t],$ for any $z\in\mathbb{T}^1,$ any $y\in Y$ and any $C_i^{(t)},$
    $$f_m(y,z,\cdot)\text{ is }(\varepsilon(t),k(t))\text{-uniformly stretching on }\{(y,z)\}\times C_i^{(t)}. $$
\end{enumerate}
We now explain why the above hypothesis of \cite[Proposition 1]{Fa02}, together with the ergodicity of $S\times R_{\alpha,\alpha'},$ implies that $\{T^t\}$ is mixing. 

Let $V\subset Y$ with $\nu(V)>0$, let $B$ be a rectangle in $\mathbb{T}^2$, and let $Q$ be the product set in $(Y\times \mathbb{T}^2)^{\widetilde{\varphi}}$ with base $V\times B\times \{0\}$ and height $\delta<\frac{1}{2}\text{ }(Q=\cup_{t=0}^{\delta}T^t(V\times B\times \{0\}).$ For $\xi>0,$ let $B_{\xi}$ be what is left from the rectangle $B$ after we have taken from its border a narrow strip of thickness $\xi$, and let $Q_{\xi}=\cup_{t=0}^{\delta}T^t(V\times B_{\xi}\times \{0\})$. (Our $B_{\xi}$ corresponds to $R_{\eta}$ in \cite{Fa02}.)

Let $\varepsilon>0$, fix a choice of $\varepsilon'>0$ (which will later be specified in terms of $\varepsilon$) and choose $\xi>0$ sufficiently small so that $(\nu\times \lambda_2)(Q_\xi)>(1-\varepsilon')(\nu\times\lambda_2)(Q).$
Fix a choice of $t>t_0$, where the conditions on $t_0$ are given below. 
    
Since $S\times R_{\alpha,\alpha'}$ is ergodic, the flow $\{T^t\}$ is also ergodic. Thus there exists $N_0>0$ and a set $E\subset (Y\times\mathbb{T}^2)^{\widetilde{\varphi}}$ with $(\nu\times\lambda_2)^{\widetilde{\varphi}}(E)>1-\varepsilon'$ consisting of points $(y,z,w,u)$ such that for $N\ge N_0$ the frequency of $T^{-s}(y,z,w,u)$ being in $Q_{\xi}$ for $0\le s\le N$ is greater than $(1-\varepsilon')^2\delta(\nu\times \lambda_2)(V\times B).$ Let $\pi:(Y\times\mathbb{T}^2)^{\widetilde{\varphi}}\to Y\times\mathbb{T}^2$ be the projection $(y,z,w,u)\mapsto (y,z,w).$ Then $G_1:=\pi(T^{-t}(E))$ satisfies $(\nu\times\lambda_2)(G_1)>1-2\varepsilon',$ because $(\nu\times\lambda_2)^{\widetilde{\varphi}}(T^{-t}E)>1-\varepsilon'$ and $\widetilde{\varphi}\ge 1/2.$  Since each point in $G_1 \times \{0\}$ can be obtained from a point in $T^{-t}E$ via a time shift of at most $3/2$, if $N\ge N_0$ where $N_0$ is sufficiently large, we have the slightly weaker frequency estimate:
\begin{equation}\label{eq:frequency}
\frac{\int_0^N \chi_{Q_{\xi}}(T^{t-s}(y,z,w,0))\ ds}{N}
\ge (1-\varepsilon')^3\delta(\nu\times\lambda_2)(V\times B),
\end{equation}
for every $(y,z,w)\in G_1$.

It also follows from the ergodicity of $S\times R_{\alpha,\alpha'}$ that for $N_0$ sufficiently large, there exists a set $G_2\subset Y\times\mathbb{T}^2$ with $(\nu\times\lambda_2)(G_2)>1-\varepsilon'$ such that for $(y,z,w)\in \pi(T^t(G_2))$, we have 
\begin{equation}\label{discrete average}
\Big|\frac{\sum_{i=1}^{N}\widetilde{\varphi}((S\times R_{\alpha,\alpha'})^{-i}(y,z,w))}{N}-1\Big| <\varepsilon',
\end{equation}
for $N\ge N_0$. Let $G=G_1\cap G_2.$
Then $(\nu\times\lambda_2)(G)>1-3\varepsilon'.$

As we will see below, equations (\ref{eq:frequency}) and (\ref{discrete average}) will replace equations (3) and (4) in \cite{Fa02}. 

Suppose $t_0>0$ is chosen so that for $t>t_0$, we have 
$$\sup_{C_i^{(t)}\in\eta_t}|C_i^{(t)}| <\min(\xi,\varepsilon'/\Vert \varphi\Vert_{C^1})$$
and 
$$\sum_{C_i^{(t)}\in\eta_t}|C_i^{(t)}|>1-\varepsilon'.$$
Assume condition (i) (the argument for condition (ii) is similar).
Then the collection of intervals $\{y\}\times C_i^{(t)}\times \{w\}$ has total $\nu\times \lambda_2$ measure greater than $1-\varepsilon'$. Hence there is a subcollection of these intervals of total $\nu\times\lambda_2$ measure greater than $1-3\sqrt{\varepsilon'}-
\varepsilon'$ such that for each interval in the subcollection, the subset of the interval consisting of points in $G$ has $\lambda_1$ measure greater than $1-\sqrt{\varepsilon'}$ times the length of the interval. For each interval in this subcollection, we truncate it by removing subintervals of length less than $\sqrt{\varepsilon'}$ times the length of the interval from the left and right sides such that the endpoints of the truncated interval are in $G$. The collection of truncated intervals obtained in this way has total $\nu\times \lambda_2$ measure greater than $1-5\sqrt{\varepsilon'}-\varepsilon'$.
We can then apply the properties of $G$ described in the previous two paragraphs at the endpoints of the truncated intervals (instead of using unique ergodicity to assume such properties at all points).  Suppose that $I=\{y_0\}\times [z_1,z_2]\times \{w_0\}$ is such a truncated interval. The points $(y_0,z_1,w_0)$ and $(y_0,z_2,w_0)$ have the same roles as $x_1$ and $x_2$, respectively, in \cite{Fa02}.  

We now describe how to change the arguments in the proof of \cite[Lemmas 1.2, 1.3, and 1.4]{Fa02} to give the estimate needed for \cite[Lemma 1]{Fa02}, namely that  $$\lambda_1( (I \times \{0\})\cap T^{-t}Q)\ge (1-\varepsilon)|I|\delta(\nu\times \lambda_2)(V\times B).$$ \cite[Lemma 1]{Fa02} will then imply that $\{T^t\}$ is mixing. 

As in \cite{Fa02}, we define 
\begin{align*}
I_m &=\{(y_0,z,w_0)\in I: 0\le t - f_m(y_0,z,w_0)\le \tilde{\varphi}((S\times R_{\alpha,\alpha'})^m(y_0,z,w_0))\}\\
    &=\{(y_0,z,w_0)\in I: N(y _0,z,w_0,t)=m\}\\
I_{m,\delta}&=\{(y_0,z,w_0)\in I_m: 0\le t-f_m(y_0,z,w_0))\le \delta\},
\end{align*}
where $N(y,z,w,t)$ is the biggest integer $m$ such that $t-f_m(y,z,w)\ge 0$, that is, the number of vertical fibers covered by $(y,z,w,0)$ during its motion under the action of the flow until time $t$. 

According to the analog of \cite[Lemma 1.2]{Fa02}, for $t$ large enough, for any $(y,z,w)\in Y\times \mathbb{T}^2$, $N(y,z,w,t)\in [\frac{1}{2}t,2t]$. This follows immediately from the fact that $\frac{1}{2}\le \widetilde{\varphi}\le \frac{3}{2},$ and we do not need any ergodicity assumption. 

As in \cite{Fa02}, we let $N_1:=N(y_0,z_1,w_0,t)$ and we consider the case that $\Delta f_{N_1|I}:=f_{N_1}(y_0,z_1,w_0)- f_{N_1}(y_0,z_2,w_0) > 0.$ 
Before truncating the intervals, we have $(\varepsilon(t),k(t))$-uniform stretch on each interval. This uniform stretch is obtained from the criteria in \cite[Lemma 2]{Fa02}. If $\varepsilon'<1/16$, then these criteria hold for the truncated intervals with $k(t)$ replaced by $\frac{1}{2}k(t),$ with no change in the $\varepsilon(t).$ That is, we obtain    $(\varepsilon(t),\frac{1}{2}k(t))$-uniform stretch on the truncated intervals.

Let $t_0>0$ be sufficiently large so that for $t>t_0,$ we have $\sup \{|C_i(t)|:C_i(t)\in\eta_t\}<\xi$, $\varepsilon(t)<1/3$ and $\frac{1}{2}k(t)\ge 4N_0.$ Then $\Delta f_{N_1|I}\ge 4N_0.$
Define 
\begin{align*}
K(I)&=\{j\in\mathbb{N}:f_j((S\times R_{\alpha,\alpha'})^{N_1}(y_0,z_2,w_0))\le \Delta f_{N_1|I}-\delta\},\\
M(I)&=\text{max }K(I).
\end{align*}
Note that here $K(I)$ is defined using just the right endpoint of $I$ instead of taking a supremum over points in $I$.

Since $|f_{j+1}-f_j|\le 3/2
$, we obtain the preliminary estimate $M(I)\ge 2N_0,$ as in the proof of \cite[Lemma 1.3]{Fa02}. 
Therefore, according to (\ref{discrete average}), the average value of $\widetilde{\varphi}$ along the backwards orbit of $(\pi(T^t(y_0,z_2,w_0))$ under $S\times R_{\alpha,\alpha'}$ for $M(I)$ iterates is approximately 1, as in the last inequality on page 450 of \cite{Fa02}. Then the analog of Lemma 1.3 stating that 
${M(I)}/{\Delta f_{{N_1}|I}}$ is approximately 1 holds. 

The analog of Lemma 1.4 is proved the same way as in \cite{Fa02}. That is, we use the fact that $M(I)$ is approximately $\Delta f_{{N_1}|I}$ to show that for $1\le j\le M(I)$, $\Delta f_{N_1+j|I}/\Delta f_{N_1|I}$ is approximately 1, provided that the length of $I$ is at most $\varepsilon'/\Vert \varphi \Vert_{C^1},$
as we are assuming. When the uniform stretch hypothesis is applied to $f_{N_1+j},$ the estimates are essentially the same as those for $f_{N_1}.$

\begin{figure}
    \includegraphics[width=\textwidth]{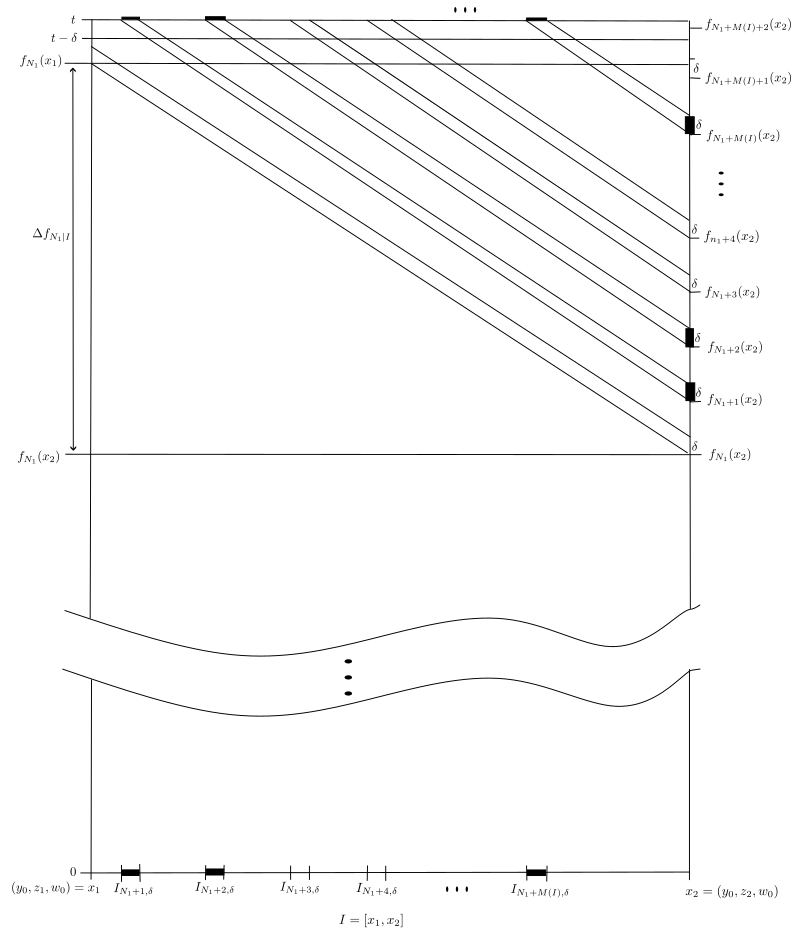}
    \captionsetup{justification=centerfirst}
        \caption{Estimating $\lambda(I \cap T^{-t}Q)$. \\
         The bold vertical segments represent the times at which the orbit of $x_2$ is assumed to be in $V\times B_{\xi}\times [0,\delta]$ for this example. The bold horizontal segments at height $t$ indicate parts of $T^t(I)$ that are in $V\times B\times [0,\delta]$. The graphs of $f_{N_1+j}$ and $f_{N_1+j}+\delta$, for $j=0,\dots,M(I),$ are approximately straight lines due to the uniform stretch condition. These graphs show times at which the orbits of points in $I$ are at height $0$ or $\delta,$ respectively, in $(Y\times\mathbb{T}^2)^{\widetilde{\varphi}}.$ The width of $I$ is actually much smaller relative to $t$ than depicted in the figure, thus making the absolute values of the slopes of the graphs of $f_{N_1+j}$ much larger.}\label{fig:fig1}
\end{figure}

The remaining part of the argument to show that the hypothesis of Lemma 1 in \cite{Fa02} holds is illustrated in Figure~\ref{fig:fig1}. Here we denote $(y_0,z_i,w_0)$ by $x_i$ for $i=1,2$. By (\ref{eq:frequency}), for $N_0$ sufficiently large, the total length of those time intervals $[f_{N_1+j}(x_2),f_{N_1+j}(x_2)+\delta]$, for $1\le j\le M(I)$, during which the orbit of $(x_2,0)$ is  in $V\times B_{\xi}\times [0,\delta]$ is at least $(1-\varepsilon')^4\delta(\nu\times \lambda_2)(V\times B)\Delta f_{N_1|I}$. For each vertical interval $[f_{N_1+j}(x_2),f_{N_1+j}(x_2)+\delta]$ over $x_2$ in Figure~\ref{fig:fig1}, the corresponding horizontal interval $I_{N_1+j,\delta}$ has length approximately equal to $\delta|I|/\Delta f_{N_1|I}.$ Thus an approximate lower bound for the total length of the horizontal intervals   $I_{N_1+j,\delta}$ such that  $T^t(I_{N_1+j,\delta}\times \{0\})\subset V\times B\times [0,\delta]$
is $(1-\varepsilon')^4\delta|I|(\nu\times\lambda_2)(V\times B_{\xi}).$ See \cite{Fa02} for the precise estimates. This shows that the hypothesis of Lemma 1 in \cite{Fa02} holds. 
Therefore, by \cite[Lemma 1]{Fa02}, $\{T^t\}$ is mixing. \qedhere
\end{proof}

As we will see in Proposition \ref{prop:ue}, the parameters in the construction of $F_1$ and $F_2$ can be chosen so that for every $\mathcal{T}\in \mathcal{T} rees$, $F_i(\mathcal{T})\times R_{\alpha,\alpha'}$ is ergodic for $i=1,2,$ and we can conclude immediately from Proposition \ref{prop:Fayad generalization} with $S=F_i(\mathcal{T})$ that 
$(F_i(\mathcal{T})\times R_{\alpha,\alpha'},\widetilde{\varphi})$ is mixing for $i=1,2.$
The proofs of both cases of Proposition \ref{prop:ue} rely on the fact that in the construction of $F_i(\mathcal{T})$, for $i=1,2$, each $(n+1)$-block consists of large numbers of consecutive appearances of each $n$-block.

\begin{prop}
\label{prop:ue}Let $F_1:\mathcal{T} rees\to\mathcal{E}$ and $F_2:\mathcal{T} rees\to \mathcal{E}\cap\text{Diff}^{\infty}(\mathbb{T}^2,\lambda_2)$ be the maps described in Sect. \ref{subsec:Review}, which satisfy Theorems \ref{AutoErgodic} and \ref{SmoothErgodic}, respectively. In the case of $F_1$, $\mathcal{E}$ denotes ergodic automorphisms of $([0,1],\lambda_1)$, and in the case of $F_2$, $\mathcal{E}\cap\text{Diff}^{\infty}(\mathbb{T}^2,\lambda_2)$ denotes $C^{\infty}$ diffeomorphisms of $\mathbb{T}^2$ that are also ergodic automorphisms of $(\mathbb{T}^2,\lambda_2). $  Then the parameters in the constructions of the maps $F_1$ and $F_2$ can be chosen so that, for every $\mathcal{T}\in\mathcal{T} rees$,  $F_1(\mathcal{T})\times R_{\alpha,\alpha'}$ and $(F_1(\mathcal{T}))^{-1}\times R_{\alpha,\alpha'}$ are ergodic with respect to $\lambda_1\times\lambda_2$, and $F_2(\mathcal{T})\times R_{\alpha,\alpha'}$ and $(F_2(\mathcal{T}))^{-1}\times R_{\alpha,\alpha'}$ are ergodic with respect to
$\lambda_2\times\lambda_2.$
 Consequently, if $\widetilde{\varphi}$ is as above, then the flows $(F_i(\mathcal{T})\times R_{\alpha,\alpha'},\widetilde{\varphi})$ and $ ((F_i(\mathcal{T}))^{-1}\times R_{\alpha,\alpha'},\widetilde{\varphi})$ are mixing on $([0,1]\times\mathbb{T}^2)^{\widetilde{\varphi}}$ with respect to $(\lambda_1\times\lambda_2)^{\widetilde{\varphi}}$ in case $i=1$ and mixing on $(\mathbb{T}^2\times \mathbb{T}^2)^{\widetilde{\varphi}}$ with respect to $(\lambda_2\times\lambda_2)^{\widetilde{\varphi}}$ in case $i=2$.
\end{prop}

\begin{proof} We will prove that $(F_i(\mathcal{T})\times R_{\alpha,\alpha'},\widetilde{\varphi})$ is ergodic. The proof that $((F_i(\mathcal{T}))^{-1}\times R_{\alpha,\alpha'},\widetilde{\varphi})$ is ergodic is similar. We will give the proof just in the case $i=2$. The case $i=1$ is the easier case, because there are no spacers if $i=1.$ In fact, if the argument below is applied in the case $i=1$, then we can show that $F_1(\mathcal{T})\times R_{\alpha,\alpha'}$ is uniquely ergodic and it suffices to apply the orginal form of \cite[Theorem 3]{Fa02}.

{\it{Case}} $i=2$. For a given $\mathcal{T}\in\mathcal{T} rees$, let $(\mathcal{S},sh,\mu)$ be the isomorphic copy of $F_2(\mathcal{T})$ described in Sect. \ref{subsec:Review}. Here $\mathcal{S}$ is a closed shift-invariant subset of $(\Sigma\cup\{b,e\})^{\mathbb{Z}}$, where $\Sigma$ is a finite alphabet not containing the ``spacer'' symbols $b$ and $e$,  $sh$ denotes the left shift, and $\mu$ is the $sh$-invariant measure obtained from the cutting and stacking model for $F_2(\mathcal{T})$. 

Suppose that for some $n\in\mathbb{N}$, we are given the collection of $n$-words $\mathcal{W}_n^c$, and suppose that each $n$-word is repeated consecutively in multiples of $\ell_n$ times, without any spacers that occur outside of $n$-words in between, every time it appears within an $(n+1)$-block. Furthermore, suppose that we are given the collection of $(n+1)$-words, $\mathcal{W}^c_{n+1}.$ Let $V$ be the set of points in $\mathcal{S}$ that have a particular position within a particular $n$-word in position 0, and let $B$ be a rectangle in $\mathbb{T}^2$. For $n\in\mathbb{N}$, let $\delta_n$ be as defined in Subsect.\ref{subsec:Review}. It suffices to show that given $\varepsilon>0$ there exist $L,M\in\mathbb{N}$ and $\delta>0$ such that if $\ell_n>L$ and $\sum_{k=n}^{\infty}\delta_k<\delta,$  then there is a set of points $x$ in $\mathcal{S}\times \mathbb{T}^2$ of measure greater than $1-\varepsilon$ such that
\begin{equation}\label{Est1}
\frac{1}{m}\sum_{k=0}^{m-1}\chi_{V\times B}\big((sh\times R_{\alpha,\alpha'})^k(x)\big)>(1-\varepsilon)\mu(V)\lambda_2(B),
\end{equation}
for all $m\ge M$.
Let $V_1$ be the set of points in $V$ such that the first $n$-word in a group of $\ell_n$ such $n$-words within an $(n+1)$-word appears in position 0. Then $\mu(V_1)=(1/\ell_n)\mu(V).$
Since each $n$-word appears the same number $f_n$ times in each $(n+1)$-word, if we choose $M$ large compared to the length of an $(n+1)$-word and we choose $\delta$ so that the measure of the set of points in $\mathcal{S}$ which have a spacer outside of the $n$-words in position 0 is small, we can arrange for 
\begin{equation}\label{Est2}
\frac{1}{m}\sum_{k=0}^{m-1}\chi_{V_1}\big((sh)^k(y)\big)>\Big(1-\frac{\varepsilon}{2}
\Big)\mu(V_1)
\end{equation}
for $y$ in a set  $A$ with $\mu(A)>1-\varepsilon$. If $(sh)^k(y)\in V_1$, and $h_n$ is the length of an $n$-word, then $(sh)^k(y),(sh)^{k+h_n}(y),(sh)^{k+2h_n}(y), \dots, (sh)^{k+(\ell_n-1)h_n}(y)$ are in $V$. Moreover, by the unique ergodicity of $R_{\alpha,\alpha'}^{h_n}$, if $\ell_n$ is sufficiently large, then for every point $(z,w)\in \mathbb{T}^2$, 
\begin{equation}\label{Est3}
\frac{1}{\ell_n}\sum_{j=0}^{\ell_n-1}\chi_B\big(R_{\alpha,\alpha'}^{jh_n}(z,w)\big)>\Big(1-\frac{\varepsilon}{2}\Big)\lambda_2(B). 
\end{equation}
By (\ref{Est2}) and (\ref{Est3}) we obtain (\ref{Est1}) for $x=(y,z,w)\in A\times\mathbb{T}^2,$ where $(\mu\times \lambda_2)(A\times\mathbb{T}^2)>1-\varepsilon.$ Here we are ignoring the edge effects coming from the fact that $(sh)^k(y)$ might be in $V_1$ for some $0\le k<m$ but some of the indices $k+h_n,\dots,k+(\ell_n-1)h_n$ might be greater than $m-1.$ We can still obtain (\ref{Est1}) if we take $M$ large compared to $\ell_nh_n$. The fact that $(1-(\varepsilon/2))^2>1-\varepsilon$ takes care of the edge effects. 
\end{proof}

\begin{rem}\label{rem:U}The first claim in Proposition \ref{prop:ue}, that is, the ergodicity of the product maps, also holds for $R_{\alpha,\alpha'}$ replaced by any homeomorphism $U$ of a compact metric space such that $U^k$ is uniquely ergodic for all $k=1,2,\dots$, and the proof is the same.  
\end{rem}

\subsection{Proof of Theorems \ref{AutoMix} and \ref{SmoothMix}}\label{Main proofs}

In the proof of Proposition \ref{prop:cancel} we will make use of the notions of towers and even equivalence.

\begin{defn} Suppose $T$ is a measure-preserving automorphism of a probability space $(X,\mathcal{B},\mu),$ and $f:X\to \mathbb{Z}^+$ is in $L^1(\mu)$. Let $X^{f}=\left\{ (x,j):x\in X,1\leq j\leq f(x)\right\} $, and let $\mu^f$ be the normalized restriction to $X^f$ of the product of $\mu$ with counting measure on $\mathbb{Z}^+.$ The {\it {tower}} $(T^f,X^f)$ over $(T,X)$ with {\it {roof function}} $f$ is the $\mu^f$-preserving map defined on $X^f$ by 
\[
T^{f}(x,j)=\begin{cases}
(x,j+1) & \text{if }j<f(x),\\
(Tx,1) & \text{if }j=f(x).
\end{cases}
\]
The roof function is said to be trivial if it is $\mu$-almost everywhere equal to $1.$ The set $X\times\{1\}$, which is often identified with $X$, is said to be the base of the tower.
\end{defn}

\begin{defn} Measure-preserving automorphisms $T$ and $S$ of  probability spaces $(X,\mathcal{B},\mu)$ and $(Y,\mathcal{C},\nu)$ are said to be {\it{evenly equivalent}} \cite{ORW82} if there exist $B\in \mathcal{B}$ and $C\in\mathcal{C}$ with $\mu(B)=\nu(C)>0$ such that the first return maps $(T_B,\mu_B
)$ and $(S_C,\mu_C)$, of $T$ to $B$ and $S$ to $C$, respectively, are isomorphic, where $\mu_B$ and $\nu_C$ are the normalized restrictions of $\mu$ to $B$ and of $\nu$ to $C$.   Equivalently, there exist measurable maps $f:X\to \mathbb{Z}^+$ and $g:Y\to \mathbb{Z}^+$ such that $\int f\ d\mu=\int g\ d\nu\  <\infty$ and $(T^f,\mu^f)$ is isomorphic to $(S^g,\nu^g).$ 
\end{defn}

The following Proposition \ref{prop:cancel} gives a cancellation property of the type discussed in Question \ref{Question} in the introduction. That is, if we assume that $U$ is as in Proposition \ref{prop:cancel}, then according to the contrapositive of the last statement in Proposition \ref{prop:cancel}, we obtain the following: If $F_i(\mathcal{T})\times U$ is isomorphic to $(F_i(\mathcal{T}))^{-1}\times U$ then $\mathcal{T}$ has an infinite branch, and therefore by Theorems 
\ref{AutoErgodic} and \ref{SmoothErgodic}, $F_i(\mathcal{T})$ and $(F_i(\mathcal{T}))^{-1}$ are isomorphic. The analogous statement holds for Kakutani equivalence.

\begin{prop}
\label{prop:cancel} Suppose that $U$ is a homeomorphism of a compact metric space such that $U$ is uniquely ergodic, and $U$ is totally ergodic with respect to the invariant measure for $U$. Furthermore, assume that there is a finite generating partition $\mathcal{R}$ for $U$ such that the measure of the boundary of each set in $\mathcal{R}$ is zero. Then it is possible to choose the parameters (depending on $U$) in the construction of the maps $F_1:\mathcal{T} rees\to \mathcal{E}$ and $F_2:\mathcal{T}rees\to\mathcal{E}\cap\text{Diff  }^{\infty}(\mathbb{T}^2,\lambda_2)$ in Theorems \ref{AutoErgodic} and \ref{SmoothErgodic} so that  $F_i(\mathcal{T})\times U$ is ergodic for $i=1,2$ and  the following implications hold for $i=1,2$: If $\mathcal{T}$ has an infinite branch, then $F_i(\mathcal{T})\times U$ is isomorphic (and hence Kakutani equivalent) to $(F_i(\mathcal{T}))^{-1}\times U$. On the other hand, if $\mathcal{T}$ does not have an infinite branch, then $F_i(\mathcal{T})\times U$ is {\emph{not}} Kakutani equivalent (and hence {\emph{not}} isomorphic) to  $(F_i(\mathcal{T}))^{-1}\times U$.  
\end{prop}

\begin{proof} Note that the unique ergodicity of $U$ and the ergodicity of $U^k$ imply that $U^k$ is uniquely ergodic. Indeed, if $\mu$ is the invariant probability measure for $U$, and $\nu$ were another ergodic invariant probability measure for $U^k$, then there would be a set $A$ with $\mu(A)=1$ and $\nu(A)=0.$ But $(1/k)(\nu+U_*\nu+\cdots +U_*^{k-1}\nu)$ is an invariant probability measure for $U$ that gives measure less than 1 to $A$ and is therefore not equal to $\mu$, contradicting the unique ergodicity of $U$. Therefore the ergodicity of $F_i(\mathcal{T})\times U$ follows from Remark \ref{rem:U} and Proposition \ref{prop:ue}.

Fix a choice of $i=1,2$ and let $F=F_i.$
The first implication in the proposition follows immediately from the fact that $F(\mathcal{T})$ and $(F(\mathcal{T}))^{-1}$ are isomorphic if $\mathcal{T}$ has an infinite branch.

Now suppose that $\mathcal{T}$ does not have an infinite branch. In order to show that 
$F(\mathcal{T})\times U$ is not Kakutani equivalent to $(F(\mathcal{T}))^{-1}\times U$ we must show the following: 
\begin{itemize}
\item[(1)] $F(\mathcal{T})\times U$ is not evenly equivalent to $(F(\mathcal{T}))^{-1}\times U.$
\item[(2)] $F(\mathcal{T})\times U$ is not isomorphic to $[(F(\mathcal{T}))^{-1}\times U]^f$ for nontrivial roof function $f$.
\item[(3)]$(F(\mathcal{T}))^{-1}\times U$ is not isomorphic to $[F(\mathcal{T})\times U]^f$ for nontrivial roof function $f$.
\end{itemize}
Conditions (2) and (3) imply that the only possible equivalence between 
$F(\mathcal{T})\times U$ and $(F(\mathcal{T}))^{-1}\times U$ is an even equivalence. But then (1) shows that an even equivalence is not possible. 
Therefore conditions (1), (2), and (3) imply that $F(\mathcal{T})\times U$ is not Kakutani equivalent to $(F(\mathcal{T}))^{-1}\times U$. 

{\it Proof of (1).} The proof of (1) is similar to the proof given in Section 9.3 of \cite{GK25} that $F({\mathcal{T}})$ and $(F({\mathcal{T}}))^{-1}$ are not evenly equivalent. Both proofs rely on Lemma 70 of \cite{GK25}, which is based on Proposition 3.2 in \cite{ORW82} and states the following: If $S_1$ and $S_2$ are evenly equivalent automorphisms with finite generating partitions $\mathcal{P}_1$ and $\mathcal{P}_2$, respectively, then for every $\varepsilon>0$ there is a finite length stationary code (with length depending on $\varepsilon$) of $(S_1,\mathcal{P}_1)$-names such that on a set of measure greater than $1-\varepsilon$ of the $(S_1,\mathcal{P}_1)$-names, the coded names are within $\varepsilon$ of $(S_2,\mathcal{P}_2)$-names in $\overline{f}$. Moreover, it follows from the proof of Lemma 70 that the measure of the set of $(S_2,\mathcal{P}_2)$-names which are approximated in this way is also greater than $1-\varepsilon.$

The main step in the proof  that $F(\mathcal{T})$ and $(F(\mathcal{T}))^{-1}$ are not evenly equivalent (for the case $F=F_1$) is  Lemma 73 of \cite{GK25}, in which we showed that any name that is finitely coded from an $s$-Feldman pattern is a bounded distance apart in $\overline{f}$ from any $s$-Feldman pattern of a different type, which leads to the negation of the conclusion of Lemma 70 in \cite{GK25} for $F(\mathcal{T})$ and $(F(\mathcal{T}))^{-1}$. The $j$th Feldman pattern is a string of symbols, $B_j$, as given in (\ref{eq:Feldman}). An $s$-Feldman pattern type refers to a particular sequence of Feldman patterns $(j_1,j_2,\dots,j_u)$ that occurs in stage $s$ of the construction in \cite{GK25}. Two such $s$-Feldman patterns $(j_1,j_2,\dots,j_u)$ and $(j_1',j_2',\dots,j_u')$ that are of different $s$-pattern type either have $\{j_1,j_2,\dots,j_u\}\cap \{j_1',j_2',\dots,j_u'\}=\emptyset$ or $(j_1,j_2,\dots,j_u)=(j_u,j_{u-1},\dots,j_1)$. (See Remark 46 in \cite{GK25}.)  In either case, a code from an $s$-Feldman pattern to an $s$-Feldman pattern of a different type (for large $u$) involves coding from one Feldman pattern (as in (\ref{eq:Feldman})) to another most of the time.

Lemma 73 is generalized to circular systems in Lemma 93 of \cite{GK25}. (This is used in the case $F=F_2$.) We will indicate the changes that need to be made to prove the analogous results for $F(\mathcal{T})\times U$ and $(F(\mathcal{T}))^{-1}\times U$. We will assume that $F=F_2$, which is the more difficult case. 

Let $\mathcal{R}$ be a finite generating partition for $U$ such that the measures of the boundaries of the sets in $\mathcal{R}$ are equal to zero. 
Each word name corresponding to $F(\mathcal{T})\times U$
is a double string of names, one for the $\mathcal{P}$ name in the $F(\mathcal{T})$ factor and 
one for the $\mathcal{R}$ name in the $U$ factor. Note that $\overline{f}$-matches between double strings of names require simultaneous matches of both symbols in some position $i$ in the first double string with the corresponding pair of symbols in some position $j$ in the second double string.

Let $\mathcal{W}^c_n$ be the collection of circular $n$-words for $F(\mathcal{T})$ as described in section 10 of \cite{GK25}, and let $q_n$ denote the length of the circular $n$-words. Let $\mathcal{V}_n$ be the collection of all $U$-$\mathcal{R}$-names of length $q_n$ that have positive probability. The collection of $n$-words for $F(\mathcal{T})\times U$ is defined to be $\mathcal{W}^c_n\times \mathcal{V}_n$. By the unique ergodicity of $U^{q_n}$ and the condition on the boundaries of the sets in $\mathcal{R}$, for $\ell_n$ sufficiently large, the frequency of each word in $\mathcal{V}_n$ in positions $ a,a+q_n,a+2q_n,\dots,a+(\ell_n-1)q_n$ of the $U$-$\mathcal{R}$-name of $y$ is approximately equal to its probability for every $y\in Y $ and every $a\in \mathbb{Z}$. Note that by the construction of circular words described in \cite[Sect.10]{GK25}, $\ell_n$ can be chosen arbitrarily large, after the words in $\mathcal{W}^c_n$ are determined.

According to \cite[Lemma 47]{GK25} each word in $\mathcal{W}_n$ that appears in 
a particular $s$-Feldman pattern appears the same number of times in each cycle. This is also true for the words in $\mathcal{W}^c_n$ for the $s$-Feldman patterns for the circular words. In addition, we know that each $\mathcal{W}^c_n$ word that appears in a particular $s$-Feldman pattern for the circular words is repeated $\ell_n$ times in a row, without any spacers in between. We define $s$-Feldman patterns and cycles in $F(\mathcal{T})\times U$-names according to the $s$-Feldman patterns and cycles of the $F(\mathcal{T})$ part of the name. By the frequency observation in the preceding paragraph, for $\ell_n$ sufficiently large, every $\mathcal{W}^c_n\times\mathcal{V}_n$ word that appears in a particular $s$-Feldman pattern of a $F(\mathcal{T})\times U$-name appears approximately the same number of times in each such cycle. This allows us to apply the Coding Lemma 43 from \cite{GK25}, to show that the image under coding of an $s$-Feldman pattern of type $r$ is bounded apart from every $s$-Feldman pattern of type $r'$ for $r\ne r'$. 

\begin{figure}
\captionsetup{justification=centerfirst}
    \includegraphics[width=\textwidth]{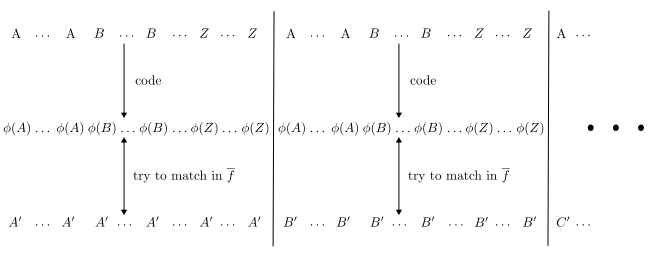}
        \caption{Failure of coding from one pattern to another.\\ If $A'\dots A'$ and $B'\dots B'$ are far apart in $\overline{f}$, then at least one of $A'\dots A'$ and $B'\dots B'$ is far apart in $\overline{f}$ from $\phi(A)\dots\phi(A)\phi(B)\dots\phi(B)\dots\phi(Z)\dots\phi(Z).$}
    \label{fig:fig2}
\end{figure}

 Figure~\ref{fig:fig2} illustrates the idea of the Coding Lemma and its generalization from $F(\mathcal{T})$-names to $F(\mathcal{T})\times U$-names. To simplify the notation we write the Feldman patterns from (\ref{eq:Feldman}) within a $F(\mathcal{T})$-name by using different capital letters for the building blocks, instead of $A$'s with subscripts. Similarly, in Feldman patterns in $(F(\mathcal{T}))^{-1}$-names, we use $A',B',C',\dots$. Any two strings with different capital letters $A',B',C',\dots$ in the third line are assumed to be a bounded distance apart in $\overline{f}$. The strings that are coded from $A,B,C,\dots$ are denoted $\phi(A),\phi(B),\phi(C),\dots$, respectively. Note that the coded string $\phi(A)$ is not completely determined by $A$ since a code $\phi$ of length $2K+1$ does not give a formula for the first $K$ and the last $K$ symbols corresponding to $A$ without considering the strings before and after $A$ (and similarly for $B,C,\dots$). However, the length of the strings $A$,$B$,$C,\dots$ can be chosen much longer than the length of the code, thereby making this edge effect negligible.

The Feldman pattern in the first and third lines of Figure~\ref{fig:fig2} come from $F(\mathcal{T})$-names and $(F(\mathcal{T}))^{-1}$-names, respectively, and the first line is assumed to have a shorter cycle than the third line. The second line shows the image of the first line under coding, and we are comparing the second and third lines. The idea is that if we happen to have a good $\overline{f}$-match between the second and the third lines for one of the constant segments (that is, repetitions of the same capital letter) in the third line, this cannot happen for any of the different constant segments, which prevents a good $\overline{f}$-match between the second and third lines. This argument would also work if the capital letters in the first string were permuted within each cycle, which is the main idea of the Coding Lemma in \cite{GK25}. 

Now we consider $F(\mathcal{T})\times U$-names and $(F(\mathcal{T}))^{-1}\times U$-names. For each $A$ in the $F(\mathcal{T})$-name, there are some number (say 3, for simplicity) possible $U$-names of the same length as $A$ that can be paired with it. We denote the paired names by $A_1$,$A_2$,$A_3$, and similarly for the other letters, $B$,$C,\dots$, as well as $A',B',C',\dots$. Then by the frequency observation from above, each of the $A_i$'s appears approximately the same number of times in each cycle (depending on $i$), and this is also true for $B_i,C_i,\dots$. Therefore the cycles in the first string are still approximately equal to permutations of each other. In the third string the $A',B',C',\dots$ also get replaced by $A_1',A_2',A_3',$ etc. Since the $A',B',C',\dots$ by themselves
are already bounded apart in $\overline{f}$, each $A_i'$ is bounded apart from every $B_j'$ by the definition of $\overline{f}$ for sequences of pairs of symbols. 
Therefore the Coding Lemma applies as before.


{\it Proof of (2).} We do not need to assume that $\mathcal{T}$ has an infinite branch to prove (2). In \cite[Sect. 9.1]{GK25} we proved that if $f$ is nontrivial, then $F(\mathcal{T})$ is not isomorphic to $(F(\mathcal{T}))^f$.
The same argument would also show that $F(\mathcal{T})$ is not isomorphic to $\big((F(\mathcal{T}))^{-1}\big)^f.$ Here we need to generalize these ideas to $F(\mathcal{T})\times U$ and $\big((F(\mathcal{T}))^{-1}\times U\big)^f$. 

Suppose that $f$ is nontrivial. If $F(\mathcal{T})\times U$ were isomorphic to $\big((F(\mathcal{T}))^{-1}\times U\big)^f$, then there would be a finite stationary code defined on  $(F(\mathcal{T})\times U,\mathcal{P}\times\mathcal{R})$-names, such that the coded names match $\big((F(\mathcal{T}))^{-1}\times U\big)^f,(\mathcal{P}\times\mathcal{R})^f\big)$-names (to be defined below) arbitrarily well in $\overline{d}$. However, we will see that this is not possible even if $\overline{d}$ is replaced by $\overline{f}.$

We define the partition $(\mathcal{P}\times \mathcal{R})^f$ to be the collection of sets in $\mathcal{P}\times \mathcal{R}$ in the base of the tower for $(\mathcal{P}\times \mathcal{R})^f$ together with the one set containing the rest of the tower, and we label points in that set by the pair of symbols $(H,H).$ Then $\big((F(\mathcal{T}))^{-1}\times U)^f,(\mathcal{P}\times \mathcal{R})^f\big)$-names consist of a double sequence of symbols corresponding to the $(\mathcal{P}\times \mathcal{R})$-names when the orbit is in the base of the tower, and $(H,H)$ for the other times. The $n$-blocks in the $\big(F(\mathcal{T})\times U,\mathcal{P}\times\mathcal{R}\big)$-names are determined by the $\big(F(\mathcal{T}),\mathcal{P}\big)$ part of the name, and the $n$-blocks for the towers start with the leftmost symbol of the $n$-blocks for the base and end right before the leftmost symbol of the next $n$-block in the base. According to \cite[Prop. 2.1, Sect. 11]{ORW82}, for most points, the length of the $n$-blocks for the towers is approximately equal to $\int f\ d(\mu\times\lambda_2)$ times the length of the $n$-blocks for the base. Since $\int f\ d(\mu\times\lambda_2)>1,$ this causes a shift of the patterns in such a way that most patterns for the base line up with a different type of pattern for the tower. 

A similar argument as in the proof of (1) shows that we cannot code from
$(F(\mathcal{T})\times U,\mathcal{P}\times\mathcal{R})$-names to $\big((F(\mathcal{T}))^{-1}\times U\big)^f,(\mathcal{P}\times\mathcal{R})^f)$-names, just using the fact that different patterns in the $(\mathcal{P}\times\mathcal{R})^f$ part of the name do not match well even in $\overline{f}$ (and we do not need to consider the $\mathcal{R}$ part of the names for this). 

Thus the argument for part (1) above shows that there is a non-zero lower bound for the $\overline{f}$ distance between such names even after coding. This proves (2). 

The proof of (3) is similar to the proof of (2).
\end{proof}

The corollary below is not needed in the proof of our main results, but we include it because it gives a new non-classifiability result for ergodic diffeomorphisms of $\T^3.$ Probably the isomorphism part could also be obtained by extending the Anosov-Katok construction used in \cite{FW22} to higher dimensions, but Corollary \ref{cor:general n} avoids this additional work. 
 
\begin{cor}\label{cor:general n}
For $n\ge 2$ and $i=1,2,$ there exist continuous maps $G_i:\mathcal{T} rees\to \mathcal{E}\cap \text{Diff}^{\infty}(\mathbb{T}^n,\lambda_n)$ such that $G_1(\mathcal{T})$ is isomorphic to $G_2(\mathcal{T})$ if and only if $\mathcal{T}$ has an infinite branch. Thus the set of pairs $(T,S)$ in $\big(\mathcal{E}\cap \text{Diff}^{\infty}(\mathbb{T}^n,\lambda_n)\big)\times \big(\mathcal{E}\cap \text{Diff}^{\infty}(\mathbb{T}^n,\lambda_n)\big)$ such that $T$ is isomorphic to $S$ is a complete analytic set, and hence not Borel. The same result is true for isomorphism replaced by Kakutani equivalence.
\end{cor}

\begin{proof}
If $n=2$, we let $G_1(\mathcal{T})=F_2(\mathcal{T})$ and $G_2(\mathcal{T})=(F_2(\mathcal{T}))^{-1}$, where $F_2$ is as in Theorem \ref{SmoothErgodic}. For $n\geq 3$, we apply Proposition \ref{prop:cancel} with $U$ equal to a rotation $R_{\alpha_1,\alpha_2,\dots,\alpha_{n-2}}$ on $\mathbb{T}^{n-2},$ where $1,\alpha_1,\dots,\alpha_{n-2}$ are linearly independent over the rationals, and let $G_1(\mathcal{T})=F_2(\mathcal{T})\times U$ and $G_2(\mathcal{T})=(F_2(\mathcal{T}))^{-1}\times U.$
    
\end{proof}

If we build a flow under $\tilde{\varphi}$ as in Proposition \ref{prop:ue}, then the identification between the top and bottom boundaries of $(\T^2\times \mathbb{T}^2)^{\tilde{\varphi}}$ depends on $F(\mathcal{T})$. By the work of S. Banerjee and P. Kunde \cite{BK19}, the map $F_2$ in Theorem \ref{SmoothErgodic} could be chosen so that the diffeomorphism $F_2(\mathcal{T})$ is real analytic, not just $C^{\infty}$, for all $\mathcal{T}\in\mathcal{T} rees$. This, together with the fact that the map $\tilde{\varphi}$ in \cite{Fa02} is real analytic, shows that the flow $(F_2(\mathcal{T}),\T^2\times \T^2,\tilde{\varphi})$ can be constructed to be a real analytic flow on a real analytic manifold.  But for the study of the complexity of classification problems, it is unsatisfactory to have the manifold depend on the transformation. To circumvent this difficulty, we apply Proposition 26 in \cite{GKpp}, to obtain diffeomorphisms of $\mathbb{T}^5$ with the usual identifications on $\mathbb{T}^5$ for all $\mathcal{T}\in\mathcal{T} rees$. It is in this transfer from $(\T^2\times \mathbb{T}^2)^{\tilde{\varphi}}$ to $\mathbb{T}^5$ that we are only able to produce $C^{\infty}$ diffeomorphisms.

\begin{prop}
\label{prop:isotopy}
For suitably chosen parameters in the construction of the map $F_2:\mathcal{T} rees\to \mathcal{E}\cap {\text{Diff}}^{\infty}(\mathbb{T}^2,\lambda_2)$ in Theorem \ref{SmoothErgodic}, there exist continuous maps $\widetilde{F}_{2,a},\widetilde{F}_{2,b}:\mathcal{T} rees\to \mathcal{ZM}\cap {\text{Diff}}^{\infty}(\mathbb{T}^5,\lambda_5)$ and a $t_0>0$ such that for each $\mathcal{T}\in\mathcal{T} rees$, $\widetilde{F}_{2,a}(\mathcal{T})$ is isomorphic to the time-$t_0$ map of the flow $(F_2(\mathcal{T})\times R_{\alpha,\alpha'},\widetilde{\varphi}),$ and $\widetilde{F}_{2,b}(\mathcal{T})$ is isomorphic to the time-$t_0$ map of the flow $((F_2(\mathcal{T}))^{-1}\times R_{\alpha,\alpha'},\widetilde{\varphi}).$
\end{prop}

\begin{proof}

We let $z=(z_1,z_2)$ and $(x,y)$ be, respectively, the coordinates on the first and second $\mathbb{T}^2$ factors in $\mathbb{T}^2\times \mathbb{T}^2$.  Let $\tilde{\varphi}(z_1,z_2,x,y)=\varphi(x,y)$ be as in Proposition \ref{prop:ue}. Let $c$ be a constant such that $\varphi(x,y)>c>0$ for all $(x,y)\in\mathbb{T}^2.$ Assume that the parameters in the construction of the map $F_2$ in Theorem \ref{SmoothErgodic} are chosen so that Propositions \ref{prop:ue} and \ref{prop:cancel} are satisfied. Also assume the map $F_2$ is chosen to satisfy the conditions imposed on $F^s$ in Proposition 26 in \cite{GKpp}. We consider the flow ($F_2(\mathcal{T})\times R_{\alpha,\alpha'},\tilde{\varphi})$, which is defined on $((\mathbb{T}^2\times\mathbb{T}^2)^{\tilde{\varphi}},\lambda_4^{\tilde{\varphi}}).$ Proposition 26 in \cite{GKpp} provides a smooth isotopy $I_u$, $0\le u\le 1$, such that $I_u\in$  Diff$^{\infty}(\mathbb{T}^2)$,  $I_u$ is the identity for $u$ in a neighborhood of $0$, and $I_u$ is $T=F(\mathcal{T})$ for $u$ in a neighborhood of 1. Moreover each $I_u$ preserves Lebesgue measure on $\mathbb{T}^2.$ We use the same isotopy here, except we
rescale it, so that it takes place on $[0,c]$ instead of $[0,1]$. Similarly to the argument in \cite{GKpp}, we can change the original flow $(T\times R_{\alpha,\alpha'},\tilde{\varphi})$ to the isomorphic flow
defined by $S^t(z,x,y,u)=(I_{u+t}(I^{-1}_u(z)),x,y,u+t)$ for $u,t\ge 0$ and $u+t\le c$; $S^t(z,x,y,u)=(z,x,y,u+t)$ for $u\ge c$, $t\ge 0$ and $u+t\le\tilde{\varphi}(x,y)$. That is, we modify the flow in Proposition \ref{prop:ue} in the first $\mathbb{T}^2$ factor for $0\le u\le c$. In the $x,y$ coordinates $(F(\mathcal{T})\times R_{\alpha,\alpha'},\tilde{\varphi})$ stays the same. The flow $S^t$ on $(\mathbb{T}^2\times\mathbb{T}^2)^{\tilde{\varphi}}$ still preserves Lebesgue measure, and the point that starts at $(z,x,y,0)$ in the base flows up to $(F(\mathcal{T})(z),x,y,\varphi(x,y))$ in the ceiling. The identification on $(\mathbb{T}^2\times\mathbb{T}^2)^{\tilde{\varphi}}$ has changed so that $(z,x,y,\varphi(x,y))$ is identified with $(z,R_{\alpha,\alpha'}(x,y),0)$. 

In \cite[Sect. 4]{Fa02}, a positive real analytic function $\phi$ on $\mathbb{T}^3$ is constructed such that the flow given by the vector field $(1/\phi)(\alpha,\alpha',1)$ takes one unit of time to go from $(x,y,0)$ to $(R_{\alpha,\alpha'}(x,y),1).$ This flow is real analytic and defined on $\mathbb{T}^3$ with the usual identifications. The flow preserves the measure $\phi\lambda$, and with this measure, it is isomorphic to $(R_{\alpha,\alpha'},\varphi)$ with measure $\lambda^{\varphi}.$ We now use the same $\phi$ on $\mathbb{T}^5,$ defined in terms of the last three coordinates, and multiply the tangent vector field for the flow $S^t$ by $1/\phi$. Then $(z,x,y,0)$ flows to $(F(\mathcal{T})(z),R_{\alpha,\alpha'}(x,y),1)$ in one unit of time and the identification on $\mathbb{T}^5$ for the reparametrized flow is the usual one. Note that $\phi$ does not depend on $F(\mathcal{T}).$ By Moser's Theorem, \cite[Theorem 5.1.27]{KH95} there is a $C^{\infty}$ diffeomorphism $\psi$ of $\mathbb{T}^5$ such that $\psi_*(\phi\lambda)=\lambda$. Define $\widetilde{F}(\mathcal{T})$ to be the composition of the time-$t_0$ map of the reparametrized flow on $\mathbb{T}^5$ with $\psi$, where $0<t_0<c$. Then $\widetilde{F}(\mathcal{T})\in$ Diff$^{\infty}(\mathbb{T}^5,\lambda).$ Moreover, $\widetilde{F}(\mathcal{T})$ is isomorphic to
the time-$t_0$ map of the flow $(F(\mathcal{T})\times R_{\alpha,\alpha'},\tilde{\varphi})$. In particular, $\widetilde{F}(\mathcal{T})\in \mathcal{ZM}.$
\end{proof}

\begin{thm} 
\label{thm:main} 
Let $\widetilde{F}_{1,a},\widetilde{F}_{1,b}:\mathcal{T} rees\to\mathcal{ZM}$ be the time $t_0$-maps of the flows $(F_1(\mathcal{T})\times R_{\alpha,\alpha'},\widetilde{\varphi})$ and $((F_1(\mathcal{T}))^{-1}\times R_{\alpha,\alpha'},\widetilde{\varphi})$, respectively, where
$F_1$ is as in Theorem \ref{AutoErgodic}, with parameters chosen to satisfy Propositions \ref{prop:ue} and \ref{prop:cancel}. Furthermore, let $\widetilde{F}_{2,a},\widetilde{F}_{2,b}:\mathcal{T} rees\to  \mathcal{ZM}\ \cap$ Diff$^{\infty}(\mathbb{T}^5,\lambda_5)$ be as in Proposition \ref{prop:isotopy}. Then the following implications hold for $i=1,2$: If $\mathcal{T}$ has an infinite branch,  then $\widetilde{F}_{i,a}(\mathcal{T})$ is isomorphic (and hence Kakutani equivalent) to $\widetilde{F}_{i,b}(\mathcal{T})$;  if $\mathcal{T}$ does not have an infinite branch, then $\widetilde{F}_{i,a}(\mathcal{T})$ is {\emph{not}} Kakutani equivalent (and hence {\emph{not}} isomorphic) to $(\widetilde{F}_{i,b}(\mathcal{T}))^{-1}.$ 
\end{thm}

\begin{proof} 
The proof is the same for $i=1$ and $i=2$. In both cases, $\widetilde{F}_{i,a}(\mathcal{T})$ and $\widetilde{F}_{i,b}(\mathcal{T})$ are isomorphic, respectively, to the time $t_0$-maps of the flows $(F_i(\mathcal{T})\times R_{\alpha,\alpha'},\widetilde{\varphi})$ and $((F_i(\mathcal{T}))^{-1}\times R_{\alpha,\alpha'},\widetilde{\varphi})$.  Thus we fix a choice of $i\in\{1,2\}$, and let $F=F_i$, $\widetilde{F}_a=\widetilde{F}_{i,a}$, and $\widetilde{F}_b=\widetilde{F}_{i,b}.$

If $\mathcal{T}$ has an infinite branch, then $F(\mathcal{T})$ is isomorphic to $(F(\mathcal{T}))^{-1}$. Since $\widetilde{\varphi}$ does not depend on the first two coordinates, $\widetilde{F}_a(\mathcal{T})$ is isomorphic to $(\widetilde{F}_b(\mathcal{T}))^{-1}.$

Now suppose that $\mathcal{T}$ does not have an infinite branch. Then by Theorem  \ref{AutoErgodic} (in case $i=1$) and Theorem \ref{SmoothErgodic} (in case $i=2$), $F(\mathcal{T})$ and $(F(\mathcal{T}))^{-1}$ are not Kakutani equivalent. 
Note that the time-$t_0$ map for the flow $(F(\mathcal{T})\times R_{\alpha,\alpha'},\tilde{\varphi})$ induced on $\mathbb{T}^4\times [0,t_0)$ is the direct product of $F(\mathcal{T})$ with a skew product having base $R_{\alpha,\alpha'}$ and a rotation by $\varphi(x,y)$ mod $t_0$ on circles of length $t_0$ in the fibers (see p. 53 in \cite{ORW82}). 
This skew product, which we will call $U$, does not depend on the first factor in the Cartesian product $F(\mathcal{T})\times R_{\alpha,\alpha'}$; in particular, it is the same for $(F(\mathcal{T}))^{-1}\times R_{\alpha,\alpha'}$. Note that $U^k$ is the skew product with base $R^k_{\alpha,\alpha'}$ and a rotation by $\beta_k(x,y):=\varphi(x,y)+\varphi(R_{\alpha,\alpha'}(x,y))+\cdots \varphi(R^{k-1}_{\alpha,\alpha'}(x,y))$ mod $t_0$ in the fibers. 

By Lemma 2.1 of \cite{Fu61}, $U$ is uniquely ergodic if it is ergodic.  Since the time-$t_0$ map of a mixing flow is ergodic, and ergodicity is preserved under inducing and taking factors, we know that $U$ is ergodic. To see that $U^k$ is ergodic, note that $U^k$ can be obtained in the same way as $U$, if we replace $\varphi$ by $\beta_k$ and we replace $R_{\alpha,\alpha'}$ by $R^k_{\alpha,\alpha'}.$ The new flow $(R_{\alpha,\alpha'}^k,\beta_k)$ on $(\mathbb{T}^2)^{\beta_k}$ still satisfies the uniform stretch condition (see Definition 1 in \cite{Fa02}), as in Proposition 1 of \cite{Fa02}. Therefore the new flow is mixing by the same argument as in \cite{Fa02}, and $U^k$ is ergodic. 

Moreover, the partition $\mathcal{R}=\{R_1\cap R_2, R_1\setminus R_2, R_2\setminus R_1, (\mathbb{T}^2\times [0,t_0))\setminus(R_1\cup R_2)\},$ where $R_1=[0,1/2)\times[0,1/2)\times[0,t_0)$ and $R_2=\mathbb{T}^2\times[0,t_0/2)$, is a generating partition for $U$ such that the measures of the boundaries of the sets in $\mathcal{R}$ are equal to zero. To see that $\mathcal{R}$ generates, note that $\{R_1, (\mathbb{T}^2\times [0,t_0))\setminus R_1\}$ generates $\mathcal{B}\times [0,t_0)$, where $\mathcal{B}$ is the collection of Borel subsets of $\mathbb{T}^2.$ Also, by the ergodicity of $U$, for every $c\in [0,t_0)$ and every $\varepsilon>0$, there are finite sequences of positive integers $(n_j)$ and Borel subsets $\{A_j\}$ of $\mathbb{T}^2$ such $\cup_{j}T^{n_j}(A_j\times [0,t_0/2))$ approximates $\mathbb{T}^2\times [c,c+t_0/2)$ (mod $t_0$) within $\varepsilon$. Thus the sets in $\mathcal{B}\times [0,t_0)$ together with $\{R_2,(\mathbb{T}^2\times [0,t_0/2))\setminus R_2\}$ generate the Borel subsets of $\mathbb{T}^2\times [0,t_0)$.

Proposition \ref{prop:cancel} and the above conditions on $U$ imply that $F(\mathcal{T})\times U$ and $(F(\mathcal{T}))^{-1}\times U$ are not Kakutani equivalent. Since $F(\mathcal{T})\times U$ and $
(F(\mathcal{T}))^{-1}\times U$ are Kakutani equivalent to $\widetilde{F}_a(\mathcal{T})$ and $\widetilde{F}_b(\mathcal{T})$, respectively, we  conclude that $\widetilde{F}_a(\mathcal{T})$ and $\widetilde{F}_b(\mathcal{T})$ are not Kakutani equivalent. 
                                                
\end{proof}

\section{\label{sec:Cancel} Cancellation of Factors from Isomorphic Cartesian Products}

As we saw in the proof of Theorem \ref{thm:main} above, for a certain family of transformations $\{F(\mathcal{T}):\mathcal{T}\in\mathcal{T} rees\}$ and a particular transformation $U$ that occurs in the proof, $F(\mathcal{T})\times U \cong F(\mathcal{T})^{-1}\times U$ if and only if $F(\mathcal{T}) \cong (F(\mathcal{T}))^{-1}$ and 
$F(\mathcal{T})\times U \sim F(\mathcal{T})^{-1}\times U$ if and only if $F(\mathcal{T}) \sim (F(\mathcal{T}
))^{-1}$. It would be interesting to know how generally this kind of cancellation holds for isomorphism and Kakutani equivalence, as stated in  Question \ref{Question} in Section \ref{sec:Introduction}.  

J.-P. Thouvenot proved \cite[Proposition 2]{T08} that for two finite entropy ergodic automorphisms $T_1$ and $T_2$ and a Bernoulli shift $B$, if $T_1\times B$ is isomorphic to $T_2\times B$, then $T_1$ is isomorphic to $T_2$. Thouvenot's result assumed that $T_1$ and $T_2$ satisfy the weak Pinsker property, but this assumption is vacuous due to the T. Austin's subsequent proof of the weak Pinsker conjecture \cite{Au18}. 

It is an open question, stated already in \cite[p. 101]{ORW82}, whether there exist automorphisms $T_1$ and $T_2$, and a finite entropy Bernoulli shift $B$ such that $T_1\times B$ and $T_2\times B$ are Kakutani equivalent, but $T_1$ and $T_2$ are not Kakutani equivalent. On the other hand, there exist finite entropy automorphisms $T_1$, $T_2$ and $U$ such that $T_1\times U$ and $T_2\times U$ are isomorphic, but $T_1$ and $T_2$ are not Kakutani equivalent. For example, we could let $T$ be a zero-entropy weakly mixing loosely Bernoulli automorphism such that $T\times T$ is not loosely Bernoulli (which exists by \cite[Chapter 14]{ORW82} \cite{R78}, \cite{R79}), and let $T_1=T$, $T_2=T\times T$, and $U=\prod_{n=1}^{\infty}T$. Then $T_1\times U$ 
is isomorphic to $T_2\times U$, but $T_1$ and $T_2$ are not Kakutani equivalent. A. Danilenko and M. Lema\'nczyk pointed out that we could also construct a similar example for any two zero-entropy transformations $T_1$ and $T_2$ that are not Kakutani equivalent, by letting $U=\prod_{n=1}^{\infty}(T_1\times T_2)$.

Now we can ask whether the analog of the above result of Thouvenot is true when $B$ is replaced by an irrational rotation. Lema\'nczyk showed that this is false. He found examples of zero entropy weakly mixing automorphisms $T_1$ and $T_2$ and an irrational rotation $R$ such that $T_1\times R$ is isomorphic to $T_2\times R$, but $T_1$ and $T_2$ are not isomorphic. In Lema\'nczyk's example, $T_1$ is Kakutani equivalent to $T_2$. We do not know if it is possible to have $T_1\times R$ isomorphic to $T_2\times R$ but $T_1$ not Kakutani equivalent to $T_2$. 

For the convenience of the reader, we give a detailed explanation of Lema\'nczyk's example. 

\medskip 
\noindent$\text{{\bf Example (due to M. Lema\'nczyk).}}$\smallskip

Let $\mathbb{T}^{1}=\{\theta:0\le\theta<1\},$ with addition mod 1,
and let $\lambda$ be Lebesgue measure on $\mathbb{T}^{1}.$ Assume
that $\alpha$ is irrational and $R_{\alpha}$ is the rotation $R_{\alpha}:\mathbb{T}^{1}\to\mathbb{T}^{1}$
given by $\theta\mapsto\theta+\alpha.$ We will show that there exist
weakly mixing transformations $T_{\varphi}$ and $T_{\psi}$ such
that $T_{\varphi}\times R_{\alpha}\cong T_{\psi}\times R_{\alpha},$
but $T_{\varphi}\ncong T_{\psi}.$ We will obtain $T_{\varphi}\times R_{\alpha}\cong T_{\psi}\times R_{\alpha}$
from Lemma~\ref{lem:IsomProd} below, and $T_{\varphi}\ncong T_{\psi}$ by Theorem~\ref{thm:NonIsom}.
The same construction also works more generally with $\mathbb{T}^{1}$
replaced by a compact abelian group.

\begin{defn}
For any measure-preserving
automorphism $T$ of a standard non-atomic probability space $(X,\mathcal{B},\mu)$
and a measurable function $h:X\to\mathbb{T}^{1}$, define the skew
product $T_{h}:X\times\mathbb{T}^{1}\to X\times\mathbb{T}^{1}$ by
$T_{h}(x,\theta)=(Tx,\theta+h(x))$. We may assume that $(X,\mathcal{B},\mu)$
is $[0,1]$ with Lebesgue measure and $\mathcal{B}=\mathcal{B}(X),$
where for any compact metric space $Y,$ we let $\mathcal{B}(Y)$
denote the Borel sets in $Y.$ 

\end{defn}

If $T$ is weakly mixing, then by \cite[Theorem 7]{JP72}
the set of functions $h$ such that $T_{h}$ is weakly mixing is a
comeager set in the space of measurable functions from $X$ to $\mathbb{T}^{1}$
with the $L^{1}(\mu)$ topology. Thus there exists a measurable $\varphi:X\to\mathbb{T}^{1}$
such that $T_{k\varphi}$ is weakly mixing for every nonzero integer
$k,$ and $T_{\psi}$ is weakly mixing, where we define $\psi=\varphi+\alpha.$

\begin{lem}\label{lem:IsomProd} Suppose $T$ is a measure-preserving
transformation of $(X,\mathcal{B}(X),\mu)$, $\alpha$ is irrational,
$\varphi:X\to\mathbb{T}^{1}$ is measurable, and $\psi=\varphi+\alpha.$
Then
\[
T_{\varphi}\times R_{\alpha}\cong T_{\psi}\times R_{\alpha}.
\]
\end{lem}

\begin{proof} Define a $(\mu\times\lambda\times\lambda$)-measure-preserving
transformation $F:X\times\mathbb{T}^{1}\times\mathbb{T}^{1}\to X\times\mathbb{T}^{1}\times\mathbb{T}^{1}$
by $F(x,u,v)=(x,u+v,v).$ Then 
\[
\begin{array}{rcl}
F((T_{\varphi}\times R_{\alpha})(x,u,v)) & = & F(Tx,u+\varphi(x),v+\alpha)\\
 & = & (Tx,u+v+\varphi(x)+\alpha,v+\alpha)\\
 & = & (T_{\psi}\times R_{\alpha})(x,u+v,v)\\
 & = & (T_{\psi}\times R_{\alpha})(F(x,u,v)).
\end{array}
\]
Therefore $F\circ(T_{\varphi}\times R_{\alpha})=(T_{\psi}\times R_{\alpha})\circ F.$ 
\end{proof}

The following lemma gives a well-known necessary condition for $T_h$ to be weakly mixing. 

\begin{lem} Suppose $T$ is a weakly mixing
measure-preserving transformation of $(X,\mathcal{B}(X),\mu),$ and
$h:X\to\mathbb{T}^{1}$ is measurable. If $T_{h}$ is weakly mixing,
then $h$ is not $T$-cohomologous to a constant, that is, there is
no measurable function $\eta:X\to\mathbb{T}^{1}$ such that $h=\eta\circ T-\eta+c,$
where $c$ is a constant. 
\end{lem}

\begin{proof}Suppose $h=\eta\circ T-\eta+c.$
Define $f:X\times\mathbb{T}^{1}\to\mathbb{C}$ by $f(x,\theta)=e^{2\pi i(-\eta(x)+\theta)}$.
Then 
\[
\begin{array}{rcl}
(f\circ T_{h})(x,\theta) & = & f(Tx,\theta+h(x))\\
 & = & e^{2\pi i(-\eta(Tx)+\theta+h(x))}\\
 & = & e^{2\pi i(-h(x)-\eta(x)+c+\theta+h(x))}\\
 & = & e^{2\pi ic}f(x,\theta).
\end{array}
\]
Therefore $f$ is a nonconstant eigenfunction of the Koopman operator
for $T_{h},$ and consequently $T_{h}$ is not weakly mixing. 
\end{proof}

\begin{lem}
\label{fourier} 
Suppose $T$ is a measure-preserving
automorphism of $(X,\mathcal{B}(X),\mu)$, $\alpha$ is irrational,
and $\varphi:X\to\mathbb{T}^{1}$ is a measurable map such that $T_{k\varphi}$
is weakly mixing for all $k\in\mathbb{Z}\setminus\{0\},$ and for
$\psi=\varphi+\alpha,$ $T_{\psi}$ is weakly mixing. Furthermore,
suppose that $\mathcal{I}:X\times\mathbb{T}^{1}\to X\times\mathbb{T}^{1}$
is an invertible measure-preserving transformation that maps vertical
fibers to vertical fibers, that is, there exists an invertible measure-preserving
transformation $S:X\to X$ such that for almost every $x\in\mathbb{T}^{1},$
we have ($\delta_{Sx}\times\lambda)\left((\mathcal{I}(\{x\}\times\mathbb{T}^{1}))\Delta(\{Sx\}\times\mathbb{T}^{1})\right)=0.$
Then $\mathcal{I}$ is not an isomorphism between $T_{\varphi}$ and
$T_{\psi}.$ 
\end{lem}

\begin{proof} Suppose that $\mathcal{I}\circ T_{\varphi}=T_{\psi}\circ\mathcal{I}.$
For $x\in X,$ let $f^{x}:\mathbb{T}^{1}\to\mathbb{T}^{1}$ be the
map such that $\mathcal{I}(x,\theta)=(Sx,f^{x}(\theta))$. It is convenient
to regard the variable $\theta$ in the domain of $f^{x}$ to be in
$\mathbb{T}^{1}$ regarded as $[0,1)$ with addition mod 1, while
regarding $f^{x}(\theta)$ to be in the isomorphic copy of $\mathbb{T}^{1}$
as the unit circle in $\mathbb{C}$ with complex multiplication. For
$\mu$-a.e. $x,$ we have a convergent Fourier series of the form
\[
f^{x}(\theta)=\sum_{n=-\infty}^{\infty}c_{n}(x)e^{2\pi in\theta},
\]
for $\lambda$-a.e. $\theta.$ Equating the second coordinates of
$(\mathcal{I}\circ T_{\varphi})(x,\theta)$ and $(T_{\psi}\circ\mathcal{I})(x,\theta)$,
where $\psi(x)=\varphi(x)+\alpha,$ we obtain $f^{Tx}(\theta+\varphi(x))=f^{x}(\theta)e^{2\pi i(\varphi(x)+\alpha)}.$
In terms of Fourier series, we have 
\[
\sum_{n=-\infty}^{\infty}c_{n}(Tx)e^{2\pi in\varphi(x)}e^{2\pi in\theta}=\sum_{n=-\infty}^{\infty}c_{n}(x)e^{2\pi i(\varphi(x))}e^{2\pi i\alpha}e^{2\pi in\theta}.
\]
Equating Fourier coefficients, we obtain 
\[
c_{n}(Tx)e^{2\pi i(n-1)\varphi(x)}=c_{n}(x)e^{2\pi i\alpha}.
\]

If we write $c_{n}(x)$ in polar form, as $c_{n}(x)=r_{n}(x)e^{2\pi i\xi_{n}(x)},$
where $r_{n}(x)\ge0,$then we have $r_{n}(Tx)=r_{n}(x).$ By ergodicity
of $T$, this implies that $r_{n}(x)$ is $\mu$-a.e. a constant $r_{n}.$
If $r_{n}\ne0,$ then we have
\begin{equation} \label{eq:rNonZero}
    e^{2\pi i\xi_{n}(Tx)}e^{2\pi i(n-1)\varphi(x)}=e^{2\pi i\xi_{n}(x)}e^{2\pi i\alpha},
\end{equation}
which implies that 
\[
\xi_{n}(Tx)-\xi_{n}(x)=-(n-1)\varphi(x)+\alpha\text{ \ (\ mod 1). }
\]
But then $(n-1)\varphi$ is cohomologous to a constant. If $n\ne1,$
this contradicts our choice of $\varphi.$ Therefore $r_{n}=0$ for
$n\ne1.$

If $n=1$ and $r_{1}\ne0,$ then equation \eqref{eq:rNonZero} implies that $e^{2\pi i\xi_{1}(x)}$
is an eigenfunction for the Koopman operator of $T,$ contradicting
the assumption that $T$ is weakly mixing. Therefore $r_{1}=0.$

Since $r_{n}=0$ for all $n,$ $f^{x}(\theta)=0$ for $(\mu\times\lambda)$-a.e.~$(x,\theta).$ This is impossible. Therefore no isomorphism of the
form $\mathcal{I}$ exists. 
\end{proof}

For Lemma \ref{lem:simple} and Theorem \ref{thm:NonIsom} we will need the following definition.

\begin{defn} A measure-preserving automorphism of $(X,\mathcal{B}(X),\mu)$ is said to be {\it{simple of order two}} if every ergodic joining of $(T,X,\mu)$ with itself is either a graph joining or an independent joining (see \cite[Chapter 12]{Gl03}).
\end{defn}

The first example of an automorphism with minimal self-joinings 
(which implies simple of order two) was constructed by D. Rudolph \cite{Ru79}. Rudolph's example is mixing of all orders. Soon afterwards, A. del Junco, M. Rahe, and L. Swanson \cite{JRS80} showed that Chacon's example of a weakly mixing but not mixing automorphism has minimal self-joinings.

\begin{lem}\label{lem:simple}
\label{vertical}
Suppose $T$ is a weakly mixing
measure-preserving automorphism of $(X,\mathcal{B}(X),\mu),$ $\alpha$
is irrational, $\varphi:X\to\mathbb{T}^{1}$ is measurable, $\psi=\varphi+\alpha,$
$T_{\varphi}$ and $T_{\psi}$ are ergodic, and $\mathcal{I}:X\times\mathbb{T}^{1}\to X\times\mathbb{T}^{1}$
is a $(\mu\times\lambda$)-measure-preserving automorphism such that
$\mathcal{I}\circ T_{\varphi}=T_{\psi}\circ\mathcal{I}.$ Assume,
in addition, that $T$ is simple of order two. Then $\mathcal{I}$
maps vertical fibers to vertical fibers (as defined in Lemma \ref{fourier}). 
\end{lem}

\begin{thm}\label{thm:NonIsom}{\bf (M. Lema\'nczyk)} Suppose $T$ is a measure-preserving
automorphism of $(X,\mathcal{B}(X),\mu),$ $\alpha$ is irrational,
$\varphi:X\to\mathbb{T}^{1}$ is measurable, and $\psi=\varphi+\alpha.$
Assume that $T$ is simple of order two, and $T_{\psi},T_{k\varphi},k\in\mathbb{Z}\setminus\{0\},$
are all weakly mixing. Then $T_{\varphi}\ncong T_{\psi}.$ 
\end{thm}

\begin{proof} It follows from Lemmas \ref{fourier} and \ref{vertical} that there
does not exist an isomorphism between $T_{\varphi}$ and $T_{\psi}.$
\end{proof}

\begin{proof}[Proof of Lemma \ref{vertical}]
Suppose\emph{ }$\mathcal{I}:X\times\mathbb{T}^{1}\to X\times\mathbb{T}^{1}$
is a $(\mu\times\lambda$)-measure-preserving automorphism such that
$\mathcal{I}\circ T_{\varphi}=T_{\psi}\circ\mathcal{I}.$ Define the
sigma-algebras $\mathcal{A}_{1}$ and $\mathcal{A}_{2}$ by 
\[
\begin{array}{ccc}
\mathcal{A}_{1} & := & \{A\times\mathbb{T}^{1}:A\in\mathcal{B}(X)\},\text{\ and}\\
\mathcal{A}_{2} & := & \{\mathcal{I}(B\times\mathbb{T}^{1}):B\in\mathcal{B}(X)\}.
\end{array}
\]
It suffices to show that $\mathcal{A}_{1}=\mathcal{A}_{2}$ (modulo
sets of $(\mu\times\lambda)$-measure zero). 

Let $\nu$ be the measure on $(X\times X,\mathcal{B}(X\times X))$
such that, for $A,B\in\mathcal{B}(X),$
\[
\nu(A\times B)=(\mu\times\lambda)\left((A\times\mathbb{T}^{1})\cap\mathcal{I}(B\times\mathbb{T}^{1})\right).
\]
 Note that 
\[
\begin{array}{rcl}
\nu(T(A)\times T(B)) & = & (\mu\times\lambda)\left((T(A)\times\mathbb{T}^{1})\cap\mathcal{I}(T(B)\times\mathbb{T}^{1})\right)\\
 & = & (\mu\times\lambda)\left((T_{\psi}(A\times\mathbb{T}^{1})\cap\mathcal{I}(T_{\varphi}(B\times\mathbb{T}^{1}))\right)\\
 & = & (\mu\times\lambda)\left((T_{\psi}(A\times\mathbb{T}^{1})\cap T_{\psi}(\mathcal{I}(B\times\mathbb{T}^{1}))\right)\\
 & = & (\mu\times\lambda)\left(T_{\psi}\left((A\times\mathbb{T}^{1})\cap(\mathcal{I}(B\times\mathbb{T}^{1}))\right)\right)\\
 & = & (\mu\times\lambda)\left((A\times\mathbb{T}^{1})\cap\mathcal{I}(B\times\mathbb{T}^{1})\right)\\
 & = & \nu(A\times B).
\end{array}
\]
Therefore $T\times T$ preserves $\nu.$ Moreover, 
\[
\begin{array}{rcl}
\nu(X\times B) & = & (\mu\times\lambda)\left((X\times\mathbb{T}^{1})\cap\mathcal{I}(B\times\mathbb{T}^{1})\right)\\
 & = & (\mu\times\lambda)(\mathcal{I}(B\times\mathbb{T}^{1}))\\
 & = & (\mu\times\lambda)(B\times\mathbb{T}^{1})\\
 & = & \mu(B),
\end{array}
\]
and similarly, $\nu(A\times X)=\mu(A).$ Hence $(T\times T,$$X\times X,\nu)$
is a joining of $(T,X,\mu)$ with itself.

For Borel sets $A,B,C,D$ in $X,$ we have 
\begin{align*}
    & \nu\left((A\times B)\cap T^{n}(C\times D)\right) \\
    = & (\mu\times\lambda)\left(((A\cap T^{n}C)\times\mathbb{T}^{1})\cap\mathcal{I}(B\cap T^{n}D)\times\mathbb{T}^{1})\right)\\
  = & (\mu\times\lambda)\left((A\times\mathbb{T}^{1})\cap\mathcal{I}(B\times\mathbb{T}^{1})\cap(T^{n}C\times\mathbb{T}^{1})\cap\mathcal{I}(T^{n}D\times\mathbb{T}^{1})\right)\\
  = & (\mu\times\lambda)\left((A\times\mathbb{T}^{1})\cap\mathcal{I}(B\times\mathbb{T}^{1})\cap T_{\psi}^{n}((C\times\mathbb{T}^{1})\cap\mathcal{I}(D\times\mathbb{T}^{1}))\right).
\end{align*}
 Thus 
\begin{align*}
&\lim_{N\to\infty}\frac{1}{N}\sum_{n=0}^{N-1}\nu\left((A\times B)\cap T^{n}(C\times D)\right)\\
=&\lim_{N\to\infty}\frac{1}{N}\sum_{n=0}^{N-1}(\mu\times\lambda)\left((A\times\mathbb{T}^{1})\cap\mathcal{I}(B\times\mathbb{T}^{1})\cap T_{\psi}^{n}((C\times\mathbb{T}^{1})\cap\mathcal{I}(D\times\mathbb{T}^{1}))\right)\\
=&(\mu\times\lambda)\left((A\times\mathbb{T}^{1})\cap\mathcal{I}(B\times\mathbb{T}^{1})\right)(\mu\times\lambda)\left((C\times\mathbb{T}^{1})\cap\mathcal{I}(D\times\mathbb{T}^{1})\right)\\
=&\nu(A\times B)\nu(C\times D),
\end{align*}
where the next-to-last equality holds because $(T_{\psi},\mu\times\lambda)$
is ergodic. Therefore $(T\times T,\nu)$ is ergodic as well. Since
$T$ is simple, $\nu$ is either a graph joining or an independent
joining. 

\noindent$\text{{\bf Case 1.}}$ Suppose $\nu$ is a graph joining.
We will show that $\mathcal{A}_{1}=\mathcal{A}_{2},$ and therefore
the conclusion of the lemma holds. 

By the definition of graph joining, there exists a $\mu$-measure-preserving
automorphism $\gamma:X\to X$ such that for any $C\in\mathcal{B}(X\times X),$
$\nu(C)=\mu(\{x\in X:(x,\gamma(x))\in C\})=\mu(\{y\in X:(\gamma^{-1}(y),y)\in C\}).$ 

Let $A\in\mathcal{B}(X).$ Then 
\[
\begin{array}{rcl}
(\mu\times\lambda)\left((A\times\mathbb{T}^{1})\cap\mathcal{I}(\gamma(A)\times\mathbb{T}^{1})\right) & = & \nu(A\times\gamma(A))\\
 & = & \mu(A)\\
 & = & (\mu\times\lambda)(A\times\mathbb{T}^{1})\\
 & = & (\mu\times\lambda)(\mathcal{I}(A\times\mathbb{T}^{1}))\\
 & = & (\mu\times\lambda)(\mathcal{I}(\gamma(A)\times\mathbb{T}^{1})).
\end{array}
\]
Therefore $A\times\mathbb{T}^{1}=\mathcal{I}(\gamma(A)\times\mathbb{T}^{1})$
(modulo a set of $(\mu\times\lambda)$-measure 0), and $\mathcal{A}_{1}\subseteq\mathcal{A}_{2}.$
If $B\in\mathcal{B}(X),$ then let $A=\gamma^{-1}(B)$ and obtain
$\mathcal{I}(B\times\mathbb{T}^{1})=\gamma^{-1}(B)\times\mathbb{T}^{1}.$
Thus $\mathcal{A}_{2}\subseteq\mathcal{A}_{1}.$ Hence $\mathcal{A}_{1}=\mathcal{A}_{2}.$

\noindent $\text{{\bf Case 2.}}$ Suppose $\nu$ is an independent
joining. We will show that this leads to a contradiction, and therefore
Case 2 does not happen. 

For $A,B\in\mathcal{B}(X),$ 
\[
\begin{array}{rcl}
(\mu\times\lambda)\left((A\times\mathbb{T}^{1})\cap\mathcal{I}(B\times\mathbb{T}^{1})\right) & = & \nu(A\times B)\\
 & = & \mu(A)\mu(B)\\
 & = & \left((\mu\times\lambda)(A\times\mathbb{T}^{1})\right)\left((\mu\times\lambda)(\mathcal{I}(B\times\mathbb{T}^{1}))\right).
\end{array}
\]
Therefore the sigma-algebras $\mathcal{A}_{1}$ and $\mathcal{A}_{2}$
are independent for $\mu\times\lambda.$ Also, $\mathcal{A}_{1},\mathcal{A}_{2},$
and $\mathcal{A}_{1}\lor\mathcal{A}_{2}$ are each invariant under
$T_{\psi}.$ 

For $n\in\mathbb{Z},$ let $h_{n}:X\times\mathbb{T}^{1}\to\mathbb{C}$
by $h_{n}(x,\theta)=e^{2\pi in\theta},$ and let $g_{n}:X\to\mathbb{C}$
by $g_{n}(x)=e^{2\pi in\psi(x)}.$ Then $h_{n}\circ T_{\psi}(x,\theta)=h_{n}(Tx,\psi(x)+\theta)=g_{n}(x)h_{n}(x,\theta).$
Since $g_{n}$ is a function of $x$ alone, $h_{n}$ is what is called
a \emph{relative eigenfunction} (that is, relative to $(T,X,\mathcal{B}(X))$)
for the Koopman operator for $T_{\psi}.$ Note that $|g_{n}(x)|=1$
for all $x\in X$ and $n\in\mathbb{Z}.$ The functions $h_{n},n\in\mathbb{Z},$
are orthonormal in $L^{2}(X\times\mathbb{T}^{1},\mu\times\lambda).$
Moreover, every function $f\in L^{2}(X\times\mathbb{T}^{1},\mu\times\lambda)$
can be written in $L^{2}(X\times\mathbb{T}^{1},\mu\times\lambda)$
as 
\begin{equation} \label{eq:f}
f=\sum_{n=-\infty}^{\infty}f_{n}h_{n},    
\end{equation}
where $f_{n}(x,\theta)=f_{n}(x)$ is a function of $x$ alone, and
$f_{n}\in L^{2}(X\times\mathbb{T}^{1},\mu\times\lambda).$ (According
to Definition 9.10 in \cite[Section 9.2]{Gl03}, it follows that
$(T_{\psi},X\times\mathbb{T}^{1},\mu\times\lambda)$ has \emph{relative
discrete spectrum} over $(T,X,\mu).)$ Taking conditional expectations
with respect to $\mathcal{A}_{1}\lor\mathcal{A}_{2}$ in~\eqref{eq:f}, we see
that for every $f\in L^{2}(X\times\mathbb{T}^{1},\mu\times\lambda),$
\begin{equation}\label{eq:Exp}
    \begin{array}{rcl}
\mathbb{E}(f|\mathcal{A}_{1}\lor\mathcal{A}_{2}) & = & \sum_{n=-\infty}^{\infty}\mathbb{E}(f_{n}h_{n}|\mathcal{A}_{1}\lor\mathcal{A}_{2})\\
 & = & \sum_{n=-\infty}^{\infty}f_{n}\mathbb{E}(h_{n}|\mathcal{A}_{1}\lor\mathcal{A}_{2}),
\end{array}
\end{equation}
where the second equality is due to the fact that $f_{n}$ is measurable
with respect to $\mathcal{A}_{1}$ and therefore with respect to $\mathcal{A}_{1}\lor\mathcal{A}_{2}$. 

Let $f\in L^{2}(X\times\mathbb{T}^{1},\mu\times\lambda)$ be measurable
with respect to $\mathcal{A}_{1}\lor\mathcal{A}_{2}$ but not measurable
with respect to $\mathcal{A}_{1}.$ (For example, $f$ could be the
characteristic function of a set in $\mathcal{A}_{2}$ with measure
strictly between $0$ and 1. By independence of $\mathcal{A}_{1}$
and $\mathcal{A}_{2},$ such a characteristic function would not be
measurable with respect to $\mathcal{A}_{1}.)$ Then it follows from~\eqref{eq:Exp} that there exists $n_{0}\in\mathbb{Z}\setminus\{0\}$ such that
$\mathbb{E}(h_{n_{0}}|\mathcal{A}_{1}\lor\mathcal{A}_{2})$ is not
measurable with respect to $\mathcal{A}_{1}.$ Note that $\mathbb{E}(h_{n_{0}}|\mathcal{A}_{1}\lor\mathcal{A}_{2})\circ T_{\psi}=g_{n_{0}}\mathbb{\mathbb{E}}(h_{n_{0}}|\mathcal{A}_{1}\lor\mathcal{A}_{2}).$
(This can be verified by integrating both sides of this equation over
$C\in\mathcal{A}_{1}\lor\mathcal{A}_{2}$ and applying the change
of variables formula.) Let $\xi=\mathbb{E}(h_{n_{0}}|\mathcal{A}_{1}\lor\mathcal{A}_{2})$
and $g=g_{n_{0}}.$ Rewriting our previous equation, we have $\xi$
measurable with respect to $\mathcal{A}_{1}\lor\mathcal{A}_{2},$
but not measurable with respect to $\mathcal{A}_{1},$ such that $\xi\circ T_{\psi}=g\cdot\xi,$
where $g$ is a $\mathcal{B}(X)$-measurable function from $X\to\mathbb{C}$
and $|g|=1.$

Let $\mathcal{P}=\{(\{x\}\times\mathbb{T}^{1})\cap\mathcal{I}(\{y\}\times\mathbb{T}^{1}):x,y\in X\}.$
Then $\mathcal{P}$ is a measurable partition of $X\times\mathbb{T}^{1}$
that generates $\mathcal{A}_{1}\lor\mathcal{A}_{2}.$ (See the definitions
in \cite[5.1.3]{VO16}.) According to Rokhlin's disintegration
theorem, there is a family $\{\omega_{P}:P\in\mathcal{P}\}$ of probability
measures on $X\times\mathbb{T}^{1}$ such that $\omega_{P}(P)=1$
for $\widehat{\omega}$-a.e. $P$ and $(\mu\times\lambda)(E)=\int\omega_{P}(E)\ d\widehat{\omega}$
for every $E\in\mathcal{B}(X\times\mathbb{T}^{1}),$ where $\widehat{\omega}$
is the measure on $\mathcal{P}$ defined by $\widehat{\omega}\left\{ \bigcup\left((\{x\}\times\mathbb{T}^{1})\cap\mathcal{I}(\{y\}\times\mathbb{T}^{1})\right):(x,y)\in F\}\right\} =(\mu\times\mu)(F)$
for $F\in\mathcal{B}(X\times X).$ Let $\pi:X\times\mathbb{T}^{1}\to X\times X$
be such that $\pi((\{x\}\times\mathbb{T}^{1})\cap\mathcal{I}(\{y\}\times\mathbb{T}^{1}))=(x,y).$
Then $\pi_{*}(\mu\times\lambda)=\nu=\mu\times\mu$ and for $S=T\times T,$
$(S,X\times X,\mathcal{B}(X\times X),\mu\times\mu)$ is a factor of
$(T_{\psi},X\times\mathbb{T}^{1},\mathcal{B}(X\times\mathbb{T}^{1}),\mu\times\lambda)$
with $S\circ\pi=\pi\circ T_{\psi}.$ (This is a special case of a
general theorem giving a correspondence between invariant subsigma
algebras for a transformation and factors of the transformation, as
described in \cite[Theorem 6.5]{EW11}.) 

Since $\xi$ is measurable with respect to $\mathcal{A}_{1}\vee\mathcal{A}_{2},$
it follows that $\xi$ is $\omega_{P}$-a.e. constant on $P$ for
$\widehat{\omega}$-a.e. $P\in\mathcal{P}.$ Thus we can define $\xi_{\pi}:X\times X\to\mathbb{C}$
by $\xi_{\pi}(x,y)=\xi\left((\{x\}\times\mathbb{T}^{1})\cap\mathcal{I}(\{y\}\times\mathbb{T}^{1})\right),$
choosing for $\widehat{\omega}$-a.e. $P\in\mathcal{P}$ the constant
value taken by $\xi$ at $\omega_{P}$-a.e. every point of $P.$ Then
$\xi_{\pi}\circ(T\times T)=g\cdot\xi_{\pi}$. However, since $\xi$
is not measurable with respect to $\mathcal{A}_{1},$ it is \emph{not}
the case that $\xi_{\pi}$ is $(\delta_{x}\times\lambda)$-a.e. constant
on $\{x\}\times\mathbb{T}^{1}$ for $\mu$-a.e. $x\in X.$ Therefore
$\xi_{\pi}$ is not measurable with respect to $\mathcal{B}(X)\times X.$

Define $\eta:X\times X\times X\to\mathbb{C}$ by $\eta(x,y_{1},y_{2})=\xi_{\pi}(x,y_{1})\overline{\xi_{\pi}(x,y_{2})}.$
Then 
\[
\begin{array}{rcl}
\eta(Tx,Ty_{1},Ty_{2}) & = & \xi_{\pi}(Tx,Ty_{1})\overline{\xi_{\pi}(Tx,Ty_{2})}\\
 & = & g(x)\xi_{\pi}(x,y_{1})\overline{g(x)}\overline{\xi_{\pi}(x,y_{2})}\\
 & = & \xi_{\pi}(x,y_{1})\overline{\xi_{\pi}(x,y_{2})}\\
 & = & \eta(x,y_{1},y_{2})
\end{array}
\]
Therefore $\eta$ is a $(T\times T\times T$)-invariant function.
Since $\xi_{\pi}$ is not measurable with respect to $\mathcal{B}(X)\times X,$
$\eta$ is not constant. This implies that $T\times T\times T$ is
not ergodic for $\mu\times\mu\times\mu$, which contradicts the assumption
that $(T,\mu)$ is weakly mixing. 
\end{proof}

\subsection*{Acknowledgements}We are grateful to Mariusz Lema\'nczyk for allowing us to include his
example, and the first author thanks him for his hospitality and patient explanations
during her visit to Torun. We also thank Adam Kanigowski for helpful conversations and for providing us with the generalization of Fayad's theorem in Proposition \ref{prop:Fayad generalization}. Moreover, we thank Matthew Foreman for asking a question that led us to Corollary \ref{cor:general n}, and we thank Jean-Paul Thouvenot and Benjamin Weiss for bringing the work of Rosenthal to our attention.

\end{document}